\documentclass[reqno]{amsart}
\usepackage{epsfig,latexsym,amsfonts,amssymb,amsmath,amscd}
\usepackage{amsthm}
\usepackage{tikz-cd}
\usepackage[mathscr]{eucal}
\usepackage{mathrsfs}
\usepackage{mathbbol}
\usepackage{array}
\usetikzlibrary{matrix,arrows,decorations.pathmorphing}
\usepackage{oldgerm,units,color}

\newcommand{\nsets}{{\mathcal P}^\sharp}
\usepackage{pst-node}
\usepackage{graphicx}

\usepackage{eepic, epic}

\usepackage{ifthen}

\def\op{{\operatorname{op}}}

\def\Amol{\langle \tT\rangle}

\def\tTz{ \tT_{\zero}}

\newcommand{\lemref}[1]{Lemma~\ref{#1}}

\def\module0{module$^\dagger$}

\def\dag{\dagger}
\def\ssemiring0{$s$-semiring$^\dagger$}

\usepackage{tikz}
\usepackage{epsfig,latexsym,amsfonts,amssymb,amsmath,amscd,graphics,epic}
\usepackage{amsfonts,amssymb,amsmath,amscd,amsthm}
\usepackage[mathscr]{eucal}
\usepackage{mathrsfs}

\usepackage{mathbbol}

\usepackage{oldgerm,units}
\usepackage{wrapfig,epsfig}

\usepackage{ifthen}
\usepackage{mathbbol}

\usepackage{amsthm}
\usepackage[mathscr]{eucal}
\usepackage{mathrsfs}
\usepackage{mathbbol}
\usepackage{oldgerm,units}
\usepackage{wrapfig}

\newtheorem{theorem}{Theorem}[section]

\newtheorem{definition}[theorem]{Definition}

\newtheorem{example}[theorem]{Example}
\newtheorem{remark}[theorem]{Remark}

\newcommand{\Net}{\mathbb N}

\newcommand{\one}{\mathbb{1}}
\newcommand{\zero}{\mathbb{0}}

\newcommand{\trop}[1]{\mathcal{#1}}

\newcommand{\tG}{\trop{G}}

\newcommand{\tT}{\trop{T}}

\newcommand{\Hom}{Hom}

    \ifx\proof\undefined
    \newenvironment{proof}{
    \smallskip
    \noindent\emph{Proof.}}{\hfill\(\Box\)
    \bigskip
    } \fi

\newcommand{\ifdef}[3]{\ifthenelse{\equal{#1}{true}}{#2}{#3}}

\definecolor{lgray}{gray}{0.90}

\def\ctw{\cdot_{\operatorname{tw}}}

\def\mcA{\mathcal A}

\def\WMor{\operatorname{WMor}}

\def\mcF{\mathcal F}
\def\mcS{\mathcal S}

\newtheorem{INote}[theorem]{Important Note}

\def\({\left(}
\def\){\right)}

\def\pipe{{\underset{{\ \, }}{\mid}}}

\def\vsemifield0{$\nu$-semifield$^\dagger$}
\def\vsemiring0{$\nu$-semiring$^\dagger$}

\def\pipe1{{\underset{{1}}{\mid}}}
\def\lmod1{\mathrel  \pipe1  \joinrel \joinrel =}

\def\CFunFF1{\operatorname{CFun} (F,F)}

\def\semiring0{semiring$^{\dagger}$}
\def\Semiring0{Semiring$^{\dagger}$}
\def\Semirings0{Semirings$^{\dagger}$}
\def\semidomain0{semidomain$^{\dagger}$}
\def\semifield0{semifield$^{\dagger}$}
\def\semifields0{semifields$^{\dagger}$}
\def\vsemifields0{$\nu$-semifields$^{\dagger}$}
\def\domain0{domain$^{\dagger}$}
\def\predomain0{pre-domain$^{\dagger}$}
\def\predomains0{pre-domains$^{\dagger}$}
\def\domains0{domains$^{\dagger}$}

\def\vdomains0{$\nu$-domains$^{\dagger}$}

\def\domains0{domains$^\dagger$}

\def\mcH{\mathcal H}
\def\mcM{\mathcal M}

\def\mcN{\mathcal N}

\newcommand{\etype}[1]{\renewcommand{\labelenumi}{(#1{enumi})}}
\def\eroman{\etype{\roman}}

\def\pipe{{\underset{{\tG}}{\mid}}}

\def\lmod{\mathrel  \pipe \joinrel \joinrel =}
\def\pipe{{\underset{{\tG}}{\mid}}}

\newtheorem{thm}[theorem]{Theorem}

\newtheorem*{thm*}{Theorem}

\newtheorem{cor}[theorem]{Corollary}

\def\Hom{\operatorname{Hom}}

\def\WMorp{\operatorname{Mor}_{\operatorname{wk}}}
\def\Morp{{\operatorname{Mor}_\preceq}}
\newtheorem{lem}[theorem]{Lemma}
\newtheorem{rem}[theorem]{Remark}
\newtheorem{prop*}{Proposition}
\newtheorem{conj*}{Conjecture}

\newtheorem{prop}[theorem]{Proposition}
\newtheorem{defn}[theorem]{Definition}
\newtheorem*{examp*}{Example}
\newtheorem*{examples*}{Examples}
\newtheorem*{remark*}{Remark}
\newtheorem*{defn*}{Definition}
\newtheorem*{note*}{Note}

\def\R{\mathbb R}
\def\C{\mathbb C}

\def\la{\lambda}

\def\tT{\mathcal T}

\numberwithin{equation}{section}

\def\M0{M_{\zero}}

\def\PS{P}

\def\Cong{\Phi}

\def\semirings0{semirings$^\dagger$}

\newcommand{\nPS}[1]{\PS_{(!#1)}}
\newcommand{\nPSo}[1]{\nPS{\one}}

\begin{document}


\title[Tensor products of bimodules]
{Tensor products  of bimodules and bimodule pairs over monoids}

\author[L.~Rowen]{Louis Rowen}
\address{Department of Mathematics, Bar-Ilan University, Ramat-Gan 52900, Israel} \email{rowen@math.biu.ac.il}

\makeatletter
\@namedef{subjclassname@2020}{%
    \textup{2020} Mathematics Subject Classification}
\makeatother


\date{\today}

\thanks{The   author  was supported by the ISF grant 1994/20 and the Anshel Peffer Chair.}

\thanks{The   author is indebted to Jaiung Jun for many enlightening conversations about the categorical aspects.}

\thanks{The   author is very grateful to the referee for many helpful suggestions concerning the original version of this paper.}


%



\makeatletter
\@namedef{subjclassname@2020}{%
    \textup{2020} Mathematics Subject Classification}
\makeatother



\begin{abstract}
We modify the well-known tensor product of modules over a semiring, to provide a general theory encompassing the classical theory of tensor products over algebras, tropical structures, and hyperfields. We discuss tensor products from three viewpoints: tensor extensions,
tensor products of modules and bimodules  over monoids, and tensor products of semi-algebras. The tensor product of residue hypermodules is functorial with respect to this construction. Special attention is paid to different kinds of morphisms and the connection to the work of Nakamura and Reyes.
\end{abstract}

\maketitle

\tableofcontents


\section{Introduction}

This is part of an ongoing project to find a general  algebraic framework that is suitable for tropical mathematics, and more generally,  idempotent semirings, and F1-geometry/algebra, and hyperfields.   In the process, we must bypass negation in the classical sense, but we still aim to recover as much of classical algebra as possible. This idea originally can be found in  \cite{Dr,Ga} (also see \cite{CC}),  and was implemented
for blueprints in \cite{Lor1,Lor2}, and put in a more general context in~\cite{Row21}, and later in \cite{BB}. A recent survey with emphasis on semirings is given in \cite{Row26}.

In the search for a broad algebraic setting, a~minimalist set of axioms, for a ``pair,'' was introduced in~\cite{JMR}, and
made more precise in \cite{AGR2}, as a module $(\mcA,+)$ over a monoid $\tT$ (which often is a group), together with a distinguished \textbf{null} submodule  $\mcA_0$,  taking the place of a zero element.
 To avoid unnecessary repetition, we refer to the introductions of \cite{AGR1} and  \cite{JMR} for more background. The primary example in this paper comes from hyperrings, cf.~\cite{V1,V2} and their modules, but as we shall see, pairs   also arise as ``tropical extensions'' and are closed under standard constructions such polynomial extensions. The most basic kinds of pairs can be found in \cite[Example~2.28]{Row26}; see  \cite{AGR2} and \cite{JMR} for  other examples.
One feature of the theory is that
pairs are accessible by means of universal algebra.

When viewing pairs categorically, there is an issue, however, with the  morphisms under consideration,  especially those other than homomorphisms that often arise in hyperfields, such as ``weak morphisms'' in which $0\in S$ implies $0\in f(S)$, and ``strong morphisms'' in which $S_1 \subseteq S_2$ implies $f(S_1)\subseteq f(S_2).$
\footnote{In the context of pairs, weak morphisms are generalized to
(Definition~\ref{wm}) and ``$\subseteq$-morphisms'' strong morphisms are generalized to (Definition~\ref{mors}). }

Much of the algebraic theory of pairs has been considered, such as linear algebra in \cite{AGR2}   roots of polynomials in \cite{Row21,Row26}, and general structure theory in \cite{JMR}.

In this paper we turn to tensor products. Tensor products arise in classical algebra in several ways: In extending the base ring (``tensor extensions''), in the character theory of groups, in the theory of simple associative algebras (in the Brauer group, whose multiplication is based on the tensor product), and in representation theory via tensor categories of modules.

An early reference for general tensor products  is \cite{Ka1}.
Our objective here is to present a tensor product for pairs which provides the constructions of the previous paragraph.
When the
morphisms are homomorphisms in the classical sense, we have a concrete tensor category  in the sense of \cite{Jag,Po}, cf.~Theorem~\ref{tenf} and its corollaries; also see Remark~\ref{adjm} for the adjoint correspondence.

However,  tensor products are not so malleable when one engages with  weak morphisms and strong morphisms, which are more appropriate for hyperfield theory.
It is not even clear what we want to take the tensor product over: the hyperfield or its power set? Since the power set of a hyperfield need not be a semiring, this is a daunting   task.

     Let us observe that the classical treatment of the tensor product of a right module and a left module over a ring $C$ (e.g., see \cite[Chapter~18]{Row08}) only utilizes the fact that $C$ is a monoid. So we shall talk of tensor products  over a multiplicative monoid, called the \textbf{underlying monoid}.
A fully satisfactory tensor theory for tensor products of monoids   over a monoid is provided in \cite{CHWW}, but it seems that the presence of the extra operation of addition often gets in the way.

Our goal in this paper is to confront the general categorical issues head on, in an effort carry a theory that includes the familiar tensor product of modules, and which is compatible with the residue (quotient) hypermodule construction originated by Krasner~\cite{krasner}.

Obstacles arise in the attempt to lift morphisms which are not homomorphisms.
One could attempt to define the tensor product of weak morphisms  via simple tensors, but  this is not  well-defined.

Recently, to bypass this difficulty, Nakamura and Reyes \cite{NakR} have tackled the issue of  tensor products in full generality for hypermagmas  (in which the hyperaddition $\boxplus$ is replaced by a binary operation $*$), in an interesting paper  which also provides intriguing examples of noncommutative hypergroups. Their solution basically is to declare that non-simple tensors are vacuous; this  is needed to obtain a general categorical tensor product of hypermagmas. The downside is that  their results do not specialize to the classical theory of tensor products of algebras. Furthermore, one may lose associativity of addition, as seen in Example~\ref{NaR1}, and basic constructions such as polynomials  cannot be viewed in terms of tensor products.
When  enlarging  the underlying monoid $\tT$  to a \textbf{tensor extension}
$ \tT'  \otimes _{\tT} \mathcal M$, cf.~Theorem~\ref{m2},   we can lift these various morphisms.

In order to be applicable to hyperfields, the theory must not require distributivity of multiplication over addition, since the power sets of certain hyperfields are not distributive. But then taking polynomials ruins associativity of multiplication. Thus, as in \cite{NakR}, for a truly encompassing theory one is necessarily led to magmas, dropping associativity in multiplication.
Nevertheless, after a brief excursion to magmas, we return to modules over monoids, which provide the most decisive results.

In a nod to the theory of tensor categories,  we conclude with a piece of the Adjoint Isomorphism Theorem in Lemma~\ref{adjc}.

\subsection{Shape of the paper}$ $

We start by reviewing the basics of  magmas, also endowed with a two-sided action by a set. The basic object is a magma $\mcA$ with an action by an underlying set $\tT$ , made into a ``pair'' with the \textbf{null} Equality is replaced by   a ``surpassing relation'' denoted as $\preceq$.
Our main object in this paper is   ``hyperpairs,'' inspired by hyperfield theory (and more generally the hypermagmas of \cite{NakR}).

There are three categories, each pertain to a major class  of morphisms: homomorphisms (the morphisms in universal algebra), $\preceq$-morphisms,   and weak morphisms.

Having prepared the groundwork,
  we bring in the main construction of this paper, the tensor product.  We start with a standard module-theoretic approach  (Definition~\ref{tp1}) which is appropriate for categories involving homomorphisms. Tensor products of homomorphisms are treated  in
 Theorems~\ref{tenf} and \ref{tenfa}].

  We obtain the  $\tT $- version of tensor extensions in Theorem~\ref{wa}, which include tropical extensions.

However,  more care is needed when  dealing with  $\preceq$-morphisms or weak morphisms in conjunction with tensor products of modules, and a partial result is given in   Theorem~\ref{adj1}.

\section{Preliminaries}


\numberwithin{equation}{section}

\subsection{Underlying algebraic structures}$ $

$\Net^*$ denotes the positive natural numbers, and we set $\Net = \Net^* \cup{\zero}.$
    \label{mag1}
\begin{definition} $ $ \eroman
\begin{enumerate}
    \item
 A \textbf{magma} is a set $\mcA$ with a binary operation $*:(\mcA \cup \{\infty\}) \times (\mcA \cup \{\infty\} ) \to \mcA\cup \{\infty\} $, not necessarily associative, satisfying $\infty * b = b*\infty = \infty$ for all $b\in \mcA.$ (Thus $\infty$ is an absorbing element.)  The magma is \textbf{total} if the operation is total, i.e.,  $*:\mcA  \times \mcA  \to \mcA  $.\footnote{In most applications the magma is total, but we could utilize $\infty$ in describing tensor products. Given any total operation, one can formally adjoin the absorbing element $\infty$.}

 In this paper a magma always has a neutral element, often denoted $\iota$, i.e., $\iota * b =b *\iota = b$ for all $b\in \mcA$.\footnote{  In \cite{NakR} these magmas are called {\it unital}, with  the neutral element  designated as $e,$ but   we reserve $e$ for~\S\ref{propN1}.}

\item
A \textbf{monoid} is a magma with an
associative operation, not necessarily commutative.  An \textbf{additive semigroup} is a commutative monoid, with the operation denoted by ``$+$," whose neutral  element is denoted as~$\zero.$

     \item
 A \textbf{pre-semiring} $\mcA$ is a monoid under  two  operations, multiplication, denoted as concatenation, and addition. We shall denote multiplication by concatenation.  We will assume that $\mcA$ has  an element ~$\zero$  that is additively neutral and also   multiplicatively absorbing, and  $(\mcA,\cdot)$ is unital, with a neutral  element~$\one$.

 \item
 A \textbf{semiring} \cite{golan92} is a pre-semiring that satisfies all the properties of a ring (including associativity and distributivity of multiplication over addition), but without
negation.

\item
A semiring $\mcA$ is a \textbf{semifield} if $(\mcA\setminus \{\zero\},\cdot)$ is a group.

\end{enumerate}
 \end{definition}


\subsubsection{$(\tT_1,\tT_2)$-bimagmas}$ $

We recall some notions from \cite{AGR2}, slightly more generally, taking \cite{NakR} into account.

 \begin{definition}\label{Tmagm}
Let $(\tT, \one)$ be a  monoid with a neutral element $\one$.
\eroman
\begin{enumerate}

  \item
A  \textbf{left $\tT$-magma}    
 is a magma
$(\mathcal A,*,\iota)$,   together with a  (left)   $\tT$-action $\tT\times \mathcal A \to \mathcal A$ (also denoted  as
concatenation), for which, \text{for all} $ a , a_i\in \tT$, $b, b_i\in \mcA$,
\begin{enumerate}

  \item $\one b = b$.

    \item $a_1(a_2 b) = (a_1a_2)b$.

    \item $\iota$ is  absorbing, i.e. $a \iota  = \iota.$

  \item The action is \textbf{distributive} over $\tT$,  in the sense that
$a(b_1*b_2) = ab_1 *ab_2.$ (In particular,   $a(b_1*b_2) = \infty$ if  $ab_1 = \infty$ or $ab_2 = \infty.$)

\end{enumerate}

     We call  the elements of $\tT$ \textbf{tangible}.  The underlying monoid $\tT$ will normally be understood from the context; in case of ambiguity, we write $\tT_{\mcA}$ for emphasis.   Then we   adjoin an absorbing element $\zero_\tT$ to $\tT $, and denote $\tTz = \tT \cup\{\zero_\tT\}$, declaring $\zero_\tT b = \iota$ for each $b\in \mcA$.

  \item A left $\tT$-\textbf{submagma} of a left $\tT$-{magma} $\mcA$ is a left submagma  closed under the $\tT$-action.

 \item A   $(\tT_1,\tT_2)$-\textbf{bimagma} $\mcA$  is  a   left $\tT_1$-magma and a right $\tT_2$-magma satisfying   $(a_1b) a_2 = a_1(b a_2)$ for all $a_i\in \tT_i$ and $b\in \mcA.$ $(\tT_1,\tT_2)$ is called the
  \textbf{underlying monoids}.

   \item A   $ \tT $-\textbf{bimagma} $\mcA$ is a $(\tT,\tT)$-{bimagma}.

\begin{enumerate}
  \item  An element  $a\in \tT$   \textbf{centralizes}     $\mcA$ if $ab = ba$ for all $b\in \mcA.$

\item    $\tT$   \textbf{centralizes}       $\mcA$ if $a$   centralizes  $\mcA$  for all $a\in \tT$.
\end{enumerate}
    \end{enumerate}
\end{definition}

\begin{rem}\label{cl0}$ $\begin{enumerate}\eroman

  \item   The \textbf{classical case} is when $\mcA$ is an algebra over a commutative ring $C$,  and $\tT =C$, viewed as a multiplicative monoid.
   \item More generally, $\tT$ could be   a partial monoid, i.e., with multiplication only defined on certain pairs of elements, with changes mutatis mutandis, using Defininition~\ref{Tmagm}(i)(b) only when the right side is defined. For example, $\tT$ could be, say, a set of matrix units with multiplication defined only on compatible matrix units, see \S\ref{FR}(2).

    \item
  If a $(\tT_1,\tT_2)$-{bimagma} $\mcA$ did not already contain a neutral absorbing element $\iota,$ we could adjoin it formally by declaring  its operation on all $b\in \mcA$ by $\iota * b =b *\iota = b$, and $a\iota =\iota a =\iota$ for all $a\in \tT.$

\end{enumerate}
\end{rem}

 (So far these definitions are rather general, and could include for example the case where $\mcA$ is a monoid with neutral element $\iota$, and $\tT$ is its semigroup of monoid endomorphisms.)

\begin{definition}\label{adm1}$ $\begin{enumerate}\eroman
\item A left  $\tT$-magma $\mcA$  is   \textbf{ weakly admissible} if $\tT\subseteq \mcA$, and, when $\mcA$ is a pre-semiring, the action of Defininition~\ref{Tmagm}(i) agrees with multiplication in $\mcA$.

 \item A  \textbf{ weakly admissible}, resp.~\textbf{admissible},   right $\tT$-bimagma   is defined analogously to (i) and~(ii).

\item A  \textbf{weakly admissible} $\tT$-bimagma $\mcA$  is a left and right weakly admissible   $\tT$-magma   (identifying both copies of $\tT $ inside $\mcA$).

\end{enumerate}
\end{definition}

\begin{lem} Let $   \mcA$ be a weakly admissible $\tT$-magma. Then $\Amol$ (the  submagma   generated by $\tT$) is an admissible $\tT$-submagma of $\mcA.$
 \end{lem}
\begin{proof}   Define a filtration on $\mcA$
by $\mcA^{(1)} = \tTz$
and, inductively, $\mcA^{(k+1)} = \mcA^{(k)}*\mcA^{(k)}.$ Then $\Amol = \cup _k \mcA^{(k)},$ and  $\iota \mcA$ implies $\mcA^{(k)}\subseteq \mcA^{(k+1)}.$ We  show by induction  that $ab\in \Amol$ for $a\in \tT$ and $b\in \Amol.$
    If $a,b\in \tT$, then $ab \in \tT$ by definition.
 Next, if $b =b_1*b_2$ for $b_1,b_2\in \mcA^{(k)}$, then $ab = a(b_1*b_2) = ab_1*ab_2\in \mcA^{(k+1)} \subseteq \Amol$.
\end{proof}


 \subsection{Prepairs and pairs}$ $

 \begin{definition}\label{symsyst}
A  $(\tT_1, \tT_2)$-\textbf{prepair}  $(\mcA,\mcA_0)$ (sometimes written  $(\mcA,\mcA_0;*)$  if we want to specify the operation $*$) is  a   $(\tT_1, \tT_2)$-bimagma $(\mcA,*,\iota)$ together with a $(\tT_1, \tT_2)$-sub-bimagma $\mcA_0$  (i.e.,   $a_1 b_0  \in \mcA_0$ and  $ b_0 a_2 \in \mcA_0$  for all $a_i\in \tT_i$  and $b_0\in \mcA_0).$

A $(\tT, \tT)$-{prepair} is also called a $\tT$-\textbf{prepair}.    We   suppress $\tT$ in the notation when it is understood.

  \end{definition}

 \begin{definition}\label{symsyst1}
We modify \cite{AGR2,JMR} slightly.
\begin{enumerate}
    \eroman

   \item
A  \textbf{pair}  $(\mcA,\mcA_0)$ is  a   weakly admissible prepair with $ab =ba$ for each $a\in \tT,$ $b\in \mcA$.\footnote{In previous work \cite{AGR1, AGR2,JMR} we assumed that all pairs are proper, to dismiss the degenerate case of $\tT = \mcA_0 = \mcA.$ But this is precisely the case treated so successfully in \cite{JM}, so we permit it here.}
 \item
A    pair $(\mcA,\mcA_0)$ is said to be \textbf{proper} if $\mcA_0\cap \tT =\emptyset$.

    \end{enumerate}
\end{definition}

\begin{INote}\label{Note1} Philosophically, $\mcA_0$ takes the place of $\zero$  in classical mathematics. The reason is that since modules need not have negation (for example, $\Net$), $\zero$ has no significant role except as a place marker in linear algebra.
\end{INote}

\subsubsection{Property N
{\cite[\S3.1]{AGR1}}}\label{propN1}$ $

  The following notion provides a crucial way to compensate for the lack of negation, using  the set $\mcA_0$ of null elements instead of $\zero.$
\begin{definition}\label{propN0}$ $
We say that a  pair
$(\mcA,\mcA_0)$ satisfies
 \textbf{Property~N}  if  there is an element $\one^\dag\in\tT$ centralizing $\mcA$ satisfying $\one *\one^\dag = \one^\dag *\one\in \mcA_0$, for which,  defining $e:= \one *\one^\dag$,  if $ a + \one \in \mcA_0$ then $a + \one = e$. (But $\one^\dag$ need not be unique.) In this case, fixing   $\one^\dag$, we denote $a^\dag = a\one^\dag$,  and $a^\circ = a* a^\dag= a^\dag *a $, for all $a\in \tT.$  
\end{definition}

  Note for $a\in \tT$ that $a^\circ = a (\one *\one^\dag ) = ae \in \mcA_0.$ 

   \begin{definition}
Let $(\mcA,\mcA_0)$ be a pair satisfying   Property~N.
 \begin{enumerate}\eroman
  \item A \textbf{quasi-zero} is an element of the form $a^\circ,$
for $a\in \tT.$

  \item $\tT^\circ := \{a^\circ : a\in \tT\}\subseteq \mcA_0$.

\item $(\mcA,\mcA_0)$  is $\circ$-\textbf{distributive} if  $(b_1 * b_2)e =b_1e  * b_2e $ for all $b_i\in \mcA.$

\item  $(\mcA,\mcA_0)$  is $\mcA_0$-\textbf{idempotent} if $b*b=b$ for all $b\in \mcA_0.$

\end{enumerate}
   \end{definition}


\subsubsection{Homomorphisms  of pairs}$ $

We consider $(\tT_1,\tT_2)$-bimagmas $\mcA$   and $\mcA'$, and   functions $f:\mcA \to \mcA'$.

\begin{definition} $ $\begin{enumerate}\eroman
  \item A  \textbf{multiplicative map} is a function $f: \mcA \to \mcA'$  satisfying $f(\iota) =\iota,$  $f(a_1b) = a_1 f(b)$, and $f(ba_2) =  f(b)a_2$  for all $a_i\in \tT_i, \ b \in \mcA.$

 \item     A  \textbf{bimagma homomorphism} (also just called a \textbf{homomorphism}) $f:\mcA \to
    \mcA'$ is a multiplicative  map satisfying $f(b_1*b_2) =  f(b_1)*f(b_2),$ $\forall b_i\in \mcA.$
\end{enumerate}
\end{definition}

\begin{definition}$ $\eroman
Suppose $(\mcA,\mcA_0)$ and $(\mcA',\mcA'_0)$ are  prepairs.
\begin{enumerate}

\item
A \textbf{paired map} $f:(\mcA,\mcA_0)\to (\mcA',\mcA'_0)$  is a multiplicative map $f:\mcA \to \mcA'$ satisfying
$f(\mcA_0)\subseteq \mcA_0'.$
\item
A \textbf{paired homomorphism}
is a paired  map which is a homomorphism.

 \end{enumerate}
\end{definition}

\begin{lem}\label{preim} $ $  \begin{enumerate}\eroman
    \item  Given a $(\tT_1,\tT_2)$-prepair $(\mcA,\mcA_0)$  and
    any   homomorphism $f: \mcA \to \mcA'$ , $(\mcA',f(\mcA_0))$ is a $(\tT_1,\tT_2)$-prepair, and then $f$ is a paired homomorphism    $f : (\mcA,\mcA_0)\to (\mcA',f(\mcA_0))$.

 \item For  any  $(\tT_1,\tT_2)$-bimodule  homomorphism $f: \mcA \to \mcA'$ where   $(\mcA',\mcA'_0)$ is a $(\tT_1,\tT_2)$-prepair, there is a prepair $(\mcA,\mcA_0)$ where $\mcA_0 = \{b\in \mcA : f(b) \in \mcA_0'\},$ and then $f$ is a paired homomorphism.
 \end{enumerate}
\end{lem}
\begin{proof} (i) By definition.

(ii)     If $a\in \tT$ and $b\in \mcA_0$ then $f(ab)= af(b) \in \mcA_0',$ so $ab\in\mcA_0.$ Then $f(\mcA_0)\subseteq \mcA'_0,$ by  definition.
\end{proof}

\subsubsection{Pre-orders and surpassing relations}$ $

 \begin{definition}\label{sur1}$ $
\begin{enumerate}\eroman
  \item  A \textbf{pre-order} on a $(\tT_1,\tT_2)$-bimagma $\mathcal A$, denoted
  $\preceq$, is a set-theoretic pre-order that respects the bimagma structure, i.e.,  for all   $b,b_i \in \mcA$:
  \begin{enumerate}

      \item $b_1\preceq b_2$ implies $a_1 b_1\preceq a_1 b_2$ and $b_1a_2 \preceq b_2a_2 $ for $a_i\in \tT_i.$
         \item   $b_i \preceq b_i'$ implies $b_1 * b_2 \preceq b_1'*b_2'.$
  \end{enumerate}

 \item  A \textbf{surpassing relation} on a weakly admissible $\tT$-bimagma $\mcA$ is a pre-order satisfying the following:
\begin{enumerate}
\item
$a_1 \preceq a_2 $ for $a_1 ,a_2 \in   \tTz$ implies $a_1 =a_2.$ (In other words, surpassing restricts to equality on~$\tTz$.)

\item
$b  \preceq \iota$ for $b \in \mcA$ implies $b=\iota$.
\end{enumerate}

      \item  A \textbf{pre-order} (resp. \textbf{surpassing relation}) on a prepair $(\mathcal A,\mcA_0)$,
 is a pre-order (resp. {surpassing relation}) $\preceq$ satisfying the condition $\iota \preceq c$ for all $c\in \mcA_0$.
\end{enumerate}
 \end{definition}

 \begin{lem}  In a prepair $(\mathcal A,\mcA_0)$ with a surpassing relation,
       $b  \preceq b* c$  and $b  \preceq c*b  $ for all $b \in \mcA$ and $c\in \mcA_0$.
 \end{lem}
 \begin{proof}
     $b \preceq b$ and $\iota \preceq c,$ so $b = b*\iota \preceq b* c.$
 \end{proof}
Surpassing relations were introduced in \cite{Row21}, and  in \cite{AGR2} for pairs, for the purposes of linear algebra.  Just as the null submodule generalizes the zero element, the surpassing relation generalizes equality.

\begin{example} Suppose $\mcA$ is a $\tT$-module.
    \begin{enumerate}
        \item If $(\mcA,\mcA_0)$ is a pair, we can define the surpassing relation $\preceq_0$ by $b_1\preceq b_2$ if and only if $b_2 = b_1*c$ for some $c\in \mcA_0.$

        \item Conversely, given a surpassing relation $\preceq$ on $\mcA$, one can define the null set $\{\mcA_0 = b\in \mcA: i\preceq b\}.$
    \end{enumerate}
\end{example}


We   insert the preorder into our categories.

\begin{definition}\label{mors} \eroman
Let $\mcA, \mcA'$ be pre-ordered  $(\tT_1,\tT_2)$-bimagmas.
  \begin{enumerate}
\item  A $\preceq$-\textbf{morphism},  (analogous to ``colax morphism'' in~\cite{NakR}) is a multiplicative {map} $f:\mcA \to \mcA'$ satisfying the following conditions:
\begin{itemize}
    \item $f(b_1) \preceq  f(b_2),\quad \forall b_1 \preceq b_2\in \mcA.$
      \item $f(b_1*b_2)\preceq f(b_1)*f(b_2),\quad \forall b_1 , b_2\in \mcA.$
\end{itemize}
      $\Morp (\mcA,\mcA')$  is the set of $\preceq$-morphisms from $\mcA$ to $\mcA'.$

      \item  A $\succeq$-\textbf{morphism}, (analogous to ``lax morphism'' in \cite{NakR}) is a multiplicative {map} $f:\mcA \to \mcA'$ satisfying the following conditions:\begin{itemize}
 \item $f(b_1) \preceq  f(b_2),\quad \forall b_1 \preceq b_2\in \mcA.$
      \item $f(b_1*b_2)\succeq f(b_1)*f(b_2),\ \forall b_1 , b_2\in \mcA.$
\end{itemize}

\end{enumerate}
\end{definition}

\begin{lem}\label{preim1}
    Given an 1:1   multiplicative map $f: \mcA \to \mcA'$ where $\mcA$
is a $\tT$-bimagma and $(\mcA',\mcA'_0)$ is a pair with surpassing relation $\preceq$, the   pair $(\mcA,\mcA_0)$ of \lemref{preim}(ii) has a surpassing relation given by $b_1\preceq b_2$ when $f(b_1) \preceq f(b_2),$ and then $f$ becomes a $\preceq$-morphism.
\end{lem}
\begin{proof}
    $a_1\preceq a_2$ implies $f(a_1)\preceq f(a_2),$ so $f(a_1)=f(a_2),$ so $a_1=a_2.$
\end{proof}

  \subsection{$\tT$-modules and bimodules}$ $

 The subsequent notions are difficult to notate for bimagmas (which need not be associative under *), so we introduce them only for modules and bimodules.
\begin{definition}$ $
    \begin{enumerate}\eroman
 \item A  left   $\tT$-\textbf{module}  $\mathcal M$ is  a  left   $\tTz$-magma  where   $*$ is also commutative and associative,
 in which case we always write $+$ in place of $*$, and $\zero$ instead of $\iota$.

 \item  A   $({\tT_1} ,{\tT_2} )$-\textbf{bimodule} is a   $({\tT_1},{\tT_2})$-{bimagma}  which is a  left   $\tT_1$-{module} and a right      $\tT_2$-{module}.

 \item Definitions~\ref{adm1}  and \ref{symsyst1}  also apply to $(\tT_1,\tT_2)$-bimodules.

\item For $(\tT_1,\tT_2)$-bimodules $\mcA$ and $
    \mcA'$,  $\Hom (\mcA,
    \mcA')$ is the set of homomorphisms  $f:\mcA \to
    \mcA'$.
    \end{enumerate}
\end{definition}

\begin{lem}\label{hom1}
    $\Hom (\mcA,
    \mcA')$ is a $(\tT_1,\tT_2)$-bimodule under the operation $(f_1  +  f_2)(b) = f_1(b) + f_2(b)$  and the actions $af: b \mapsto a f(b)$ and $fa: b \mapsto  f(b)a$.
\end{lem}
\begin{proof}
Standard.
\end{proof}
\begin{remark} For $(\tT_1,\tT_2)$-bimodules,
    $\Hom (\mcA,
    \mcA')$ is a $(\tT_1,\tT_2)$-bimodule under the action of Remark~\ref{hom1}.
\end{remark}

  Just as  the category theory for (semi)algebras differs from the category theory for bimodules, we want to differentiate between $({\tT_1} ,{\tT_2} )$-{bimodules}  and $\tT$-{semirings}  (resp. pre-semirings), which we now define.

\begin{definition}$ $
    \begin{enumerate}\eroman
 \item
 A $ \tT   $-\textbf{pre-semiring} is a    $ \tT$-bimodule $(\mcA,+)$,  which  also  has a second operation which we denote as $\cdot$, with $a_1(b_1\cdot b_2)= (a_1 b_1)\cdot b_2$ and $(b_1\cdot b_2)a_2=   b_1 \cdot (b_2 a_2)$ for $a_i\in \tT,$ $b_1,b_2\in \mcA.$ In particular, when $\mcA$ is weakly admissible,  $a \cdot b = ab$ and $b\cdot a = ba$ for all $a \in \tT$ and $b\in\mcA.$

    \item A $ \tT  $-\textbf{semiring} is a $ \tT  $-pre-semiring
   which is a semiring, i.e., $\cdot$ is associative and distributes over~$+$.

  \item  A  \textbf{pre-semiring homomorphism} $f:\mcA \to
    \mcA'$  of  $ \tT  $-pre-semirings is a homomorphism    also satisfying  $f(b_1b_2) = f(b_1) f(b_2)$ for all $ b_i \in \mcA.$

\item
A pre-semiring    $\mathcal A$ with a surpassing relation $\preceq$ is $\preceq$-\textbf{distributive} if $b(b_1+b_2)\preceq b b_1  \, +\, b b_2 $ and $(b_1+b_2)b\preceq  b_1 b \, +\, b_2b $ for all $b, b_i \in \mcA.$

    \end{enumerate}
\end{definition}

 \subsubsection{Negation maps}$ $

At times we can define the negation map, the mainstay of
\cite{Row21}. A \textbf{negation map} $(-)$ on a $\tT$-module $(\mcA,\mcA_0)$ is
 an additive automorphism $(-)$  of order $\le 2$  also
 defined on $\tT$, such that $$(-)(ab) = ((-)a)b =
a((-)b),\qquad  \forall a\in \tT, \quad b\in \mcA,$$  and
$(-)\mcA_0 = \mcA_0.$

We write $b_1(-)b_2$ for $b_1+((-)b_2).$ Thus  $b^\circ = b(-)b,$ and
 $\mcA_0$ contains the set  $\mcA^\circ =
\{b^\circ: b\in \mcA\}.$
Often  $\mcA_0 = \mcA^\circ.$

We say $\mcA$ is \textbf{uniquely negated} if $\one +a =e$ implies $a =(-)\one.$

\begin{lem}
The negation map satisfies
$(-) b= ((-)\one)b.$
\end{lem}
\begin{proof}
$(-) b=(-)(\one b)=  ((-)\one)b.$
\end{proof}

Hence, perhaps surprisingly,   if $b_1\preceq b_2,$ then $(-)b_1\preceq (-)b_2$.
Hence, if $b_1\preceq b_2$   then $b_2 (-) b_1  \succeq \iota$ and
$b_1 (-) b_2 \succeq \iota,$ cf.~\cite[Lemma~2.11]{JMR}.

 \subsubsection{Weak morphisms}
\begin{definition}\label{wm}
A \textbf{weak morphism of bimodule prepairs} is a paired multiplicative map $f:(\mcA,\mcA_0)\to (\mcA',\mcA'_0)$, satisfying $\sum b_i \in \mcA_0$ for $b_i\in \mcA$ implies $\sum f(b_i)\in \mcA'_0$.

 $\WMorp (\mcA,\mcA')$ denotes the set of weak morphisms from $(\mcA,\mcA_0)$ to~$(\mcA',\mcA'_0).$
\end{definition}
\begin{lem}
For bimodule prepairs $(\mcA,\mcA_0)$ and $(\mcA',\mcA'_0),$ define $ \WMorp(\mcA,\mcA')_0$ to be the weak morphisms for which $f(\mcA)\subseteq \mcA'_0.$
\begin{enumerate}  \eroman   \item
$(\WMorp(\mcA,\mcA'),\WMorp(\mcA,\mcA')_0)$ is a bimodule prepair.

 \item For  pairs $(\mcA,\mcA_0)$ and $(\mcA',\mcA'_0),$
$(\WMorp(\mcA,\mcA'),\WMorp(\mcA,\mcA')_0)$ is a bimodule  pair, where we define $(f_1+f_2)(b) = f_1(b)+f_2(b),$  $(af)(b)= af(b)$, $(fa)(b) = f(b)a$.
\end{enumerate}
\end{lem}
\begin{proof} (i) is as in Lemma~\ref{hom1}.

    The main verification for (ii) is that the sum $f_1 + f_2$ of two weak morphisms is a weak morphism. If $\sum b_i \in \mcA_0,$ then
    $(f_1+f_2)(\sum b_i) = f_1(\sum b_i) +f_2(\sum b_i) \in \mcA_0+\mcA_0 \subseteq \mcA_0.$
\end{proof}

\begin{lem}[As in {\cite[Lemma~2.10]{AGR1}}]
    Every  $\preceq$-{morphism} of bimodule prepairs is a weak morphism.
\end{lem}

  \subsection{Major Examples}\label{ME}

  For the reader's convenience, we review examples of pairs, drawing heavily from~\cite{Row26}.

\begin{enumerate}\eroman

   \item    The prototype which we call the \textbf{classical pair} and whose structure theory we  imitate, is the pair $(\mcA,\{\zero\}),$ where $\mcA$ is a faithful algebra over an integral domain $C.$
$\tT$ is the monoid $C\setminus \{0\}$ and $\tTz=C,$ and $\mcA_0 = \{0\}.$

   \item If $P$ is a prime ideal of an integral domain $\mcA,$ we have the pair $(\mcA,P)$ with
 underlying monoid  $C\setminus P$, which is a domain that provides a different perspective to the classical integral domain $\mcA/P$.

   \item  In the same spirit as (i) and (ii), we have the pair $(\mcA,\{\zero\}),$ where $\mcA$ is any semiring, and $\tT$ is any submonoid, such as $\{1\}.$

 \item  Any semiring becomes a uniquely negated pair, when we put $\tT = \{1\}$ and $\mcA_0 = \{ b+b:b\in \mcA\}$. This has the flavor of characteristic 2.

   \item \textbf{Hyperpairs}, perhaps the most interesting nonclassical example, are pairs obtained from the power set of a hyperring, as explained below in \S\ref{hype}.

  \item  The \textbf{tropical extension} $ \mathcal L \rtimes \tG$ of a semigroup $\mathcal L$ by an ordered abelian group $(\tG,+)$
consists of the set $\mathcal L\times \tG$, made into a semigroup with the following addition:
   \begin{equation}
\label{basicex17}(\ell_1,b_1) + (\ell_2,b_2) = \begin{cases}
(\ell_1,b_1) \text{ if } b_1 > b_2,
 \\ (\ell_2,b_2) \text{ if } b_1 < b_2,  \\  (\ell_1 + \ell_2,\, b_1)
 \text{ if }   b_1 =  b_2  .
\end{cases}\end{equation}

When  $(\mathcal L,\mathcal L_0)$ is a semiring pair we make $ \mathcal L \rtimes \tG$ into a semiring   under the componentwise operation $(\ell_1,b_1) (\ell_2,b_2) = (\ell_1\ell_2,b_1 + b_2)$. Then $(\mathcal L \rtimes \tG,  \mathcal L_\zero\rtimes \tG)$ is a  semiring    pair with underlying monoid $\tT\times \tG$,  under the action $(a,g)(\ell,g') = (a\ell,g+g')$ for $a\in \tT,$ $\ell \in \mathcal L,$ $g,g'\in \tG$.   \item \textbf{Supertropical pairs}. Suppose $\tG$ is an ordered abelian  monoid with absorbing minimal element $\zero_{\tG}$, and $\tT$ is a monoid, together with a monoid homomorphism $\mu: \tT\to \tG$. The  \textbf{supertropical pair},  the foundation of supertropical algebra, \cite{IR}, is $\mcA =\tG \cup \tG^\circ \cup \{-\infty\},$ where $\tT$ is an ordered abelian group and $\mcA_0 = \tG^\circ$, with addition given by $a_1+a_2 = \max(a_1,a_2)$ for $a_1\ne a_2$ and $a+a = a^\circ.$  (Thus $-\infty^\circ = -\infty$.) This can be viewed as a special case of a supertropical extension.

One instance is through \textbf{tropicalization} of a field $F$ with a nonarchimedian valuation $v$. Then the value group $v(F)$ yields a supertropical pair $(\mcA,\mcA_0)$ where $\mcA =(v(F)\cup v(F)^\circ,$ as in the previous paragraph,
and the valuation $v$ induces a $\preceq_0$-morphism from the classical pair $(F,\{0\})$ to $(\mcA,\mcA_0)$ since if $\sum f_i = 0 $ then $v(\sum f_i) =-\infty \preceq_0 \sum v(f_i).$

\item Any $\tTz$-module $\mcA$ gives rise to the \textbf{doubled $\tTz$-module} $\widehat{\mcA} = \mcA\times \mcA$ with underlying monoid $(\tT \times \zero)\cup (\zero\times \tT),$    and    null ideal $\hat \mcA_0 = \{(b,b): b\in \mcA.$     When $\mcA$ is a semiring, define
  \textbf{Twist multiplication} on $\widehat{\mathcal A}$
 as follows:
\begin{equation}\label{twis} (b_1,b_2)\ctw (b'_1,b'_2) =
 (b_1b'_1 + b_2 b'_2, b_1 b'_2 + b_2 b'_1),\qquad (b_1,b_2), (b'_1,b'_2) \in \widehat{\mathcal A}.  \end{equation}
$(\hat \mcA,\hat \mcA_0)$ has the negation map $(-)(b_1,b_2)=(b_2,b_1)$.

   \item  Examples of pairs satisfying Property N which are not uniquely negated, taken from \cite{AGR2}:

    \begin{itemize}
   \item \cite[Example~2.21]{AGR2}  is the semiring pair $(\mcA ,\mcA_0 )$ with $\mcA = \tTz \cup \{e\}$, $\tT$~an arbitrary monoid, and $\mcA_0 = \{ \zero, e\},$
 satisfying $b_1+b_2= e$ for all $b_1,b_2 \in \tT \cup \{e\}.$ The pair $(\mcA ,\mcA_0 )$ satisfies Property N, and is \textbf{terminal} in the sense that all pairs map to it by sending $b\mapsto e$ for all $b\ne \zero$.

     \item  A terminal idempotent $\mcA_0$-bipotent pair:
 Let $\mcA = \tTz \cup \{\infty\}$ where   $a+a = a$, but $a+a' =\infty$ for each $a\ne a'\in \tT$; then $\mcA_0 = \{\zero, \infty\}.$
    \end{itemize}
   \end{enumerate}

In addition to covering a range of examples, pairs permit useful constructions of new pairs:

 Suppose $(\mcA,\mcA_0)$ is a semiring pair. \begin{itemize}
\item   For any $\tT$-module $\mcA$ and monoid $\mcM,$ we define  $\mcA[\mcM]$ to be formal sums $\sum_{u\in \mcM, \, a_u \in \tT } a_u u$,  under componentwise addition $ \sum a_u u + \sum b_u u =\sum (a_u+b_u)u$,  which  is a module over the underlying monoid $\tT_\mcM = \{ a u : a\in \tT, u \in \mcM\},$
    under the action $(au)\sum a_v v = \sum(a a_v)uv.$
    When $\mcA$ is a semiring, $\mcA[\mcM]$ becomes a semiring, called the \textbf{monoid semiring}    with the \textbf{convolution product} $(\sum a_u u)(\sum a'_v v)= ((\sum _{uv=s} a_u a'_v)s),$ provided this is defined.   The special case $(\mcA[\Lambda],\mcA_0[\Lambda])$ is called the \textbf{polynomial pair}.

 \item For an additive semigroup  $\mcA$,   $\mcA^S$ is an additive semigroup under elementwise addition of functions,  whose zero element is the function with empty support, i.e., sending   all   elements to~$\zero$. The elements of $\mcA^S$ can be viewed as $S$-tuples,
The \textbf{function pair}  $(\mcA^S,\mcA_0^S)$ is a   pair with underlying  monoid defined as the elements of $\tT^S$ having support of order 1.

\item The sub-semigroup $\mcA^{(S)}$ of functions having finite support
  is the \textbf{semigroup direct sum} and defines a sub-pair of $(\mcA^S,\mcA_0^S)$.

 \item Likewise, there is the matrix pair over a pair, defined in the natural way.
   \end{itemize}

 For example, although polynomials over hyperrings are difficult to describe and are not hyperrings, the  polynomial pair over a hyperpair is a pair. Polynomial pairs will be seen to be tensor extensions of pairs.

 \section{Hypermagmas and hyperpairs}\label{hype}$ $

 We follow the treatment of Nakamura and Reyes \cite{NakR},   casting hypergroups into a broader context which has many  more interesting examples.

\begin{defn}[\cite{NakR}, generalizing \cite{Mar}]\label{Hyp00} $ $
\begin{enumerate}\eroman
    \item $\mathcal{P}(\mcH)$ denotes the power set of a set $\mcH,$
    and $\nsets(\mcH) = \mathcal{P}(\mcH) \setminus \emptyset.$

A  \textbf{hypermagma} $(\mathcal{H},*)$ is a set with  a binary operation $*:\mcH \times \mcH \to \mathcal{P}(\mcH)$,  permitting the product of two elements to be the empty set\footnote{Customarily one requires $a_1* a_2 \in \nsets(\mcH)$ for all $a_1,a_2\in \mcH.$ But one could formally adjoin $\emptyset$, so the  Nakamura-Reyes definition subsumes the customary definition, and we shall need it in what follows. $\emptyset$ could be interpreted  as ``undefined.''}, together with a left and right $\mcH$-action whose natural extension to  $ \mathcal{P}(\mcH)$ makes $ (\mathcal{P}(\mcH),*)$ an $\mcH$-bimagma when we define, for     $S_1,S_2\in \mathcal P(\mathcal{H}),$   $$S_1 * S _2:= \cup _{s_i \in S_i}\  s_1 * s_2, \qquad \emptyset * S = S * \emptyset = \emptyset .$$ We view $\mcH \subseteq \mathcal{P}(\mcH)$ by identifying $a$ with $\{ a\}$.

From now on we follow the customary use of $\boxplus$ and $\zero$ instead of $*$ and $\iota$, to denote that the operation is associative.
We call $\boxplus$ ``hyperaddition.''

\item
  A  \textbf{hyperzero} $\zero$ (if it exists) satisfies   $\zero \boxplus a = a = a \boxplus \zero,$ for every $a\in \mathcal{H}.$

 \item A \textbf{hypersemigroup} is a hypermagma
$(\mathcal{H},\boxplus,\zero)$, where
\begin{enumerate}\eroman
   \item The hyperaddition
$\boxplus$ is  associative in the sense that   $(a_1
\boxplus a_2) \boxplus a_3 = a_1 \boxplus (a_2\boxplus a_3)$ for all
$a_i$ in $\mathcal{H}.$
  \item
$\zero \in \mathcal{H}$ is the hyperzero.
 \end{enumerate}

 \item
 A \textbf{hypernegative} of an element $a$ in   $(\mathcal{H},\boxplus,\zero)$ (if it exists) is an
element  ``$-a$'' for which $\zero \in a \boxplus (-a)$ and $\zero \in (-a) \boxplus a$.

 \item If the hypernegative $-\one$ exists in $\mcH,$ then we define $e = \one \boxplus (-\one).$

  \item
A \textbf{hypergroup} is a   hypersemigroup $(\mathcal{H},\boxplus,\zero)$ for which every element~$a$
has a unique \textbf{hypernegative} denoted $-a$, whereby, for all $a_i\in \mcH,$

\begin{enumerate}
   \item $(-)(a_1\boxplus a_2) = (-)a_2 \boxplus (-)a_1.$
   \item $-(-a_1) =a_1.$
  \item $\mcH$ is
\textbf{reversible} in the following sense:

$ a_3 \in a_1 \boxplus a_2$    iff  $ a_2 \in a_3
\boxplus (-a_1)$.\footnote{In \cite{NakR}, a \textbf{mosaic} is a reversible  with a hyperzero.}
 \end{enumerate}

 \item
A \textbf{hypersemiring} (resp.~\textbf{hyperring}) is an  abelian\footnote{In \cite{NakR} this is called ``canonical.''} hypersemigroup (resp.~hypergroup) $(\mathcal{H},\boxplus,\zero)$, where $\mcH$ has multiplication distributing over hyperaddition, providing  $\mathcal P (\mathcal H  )$ with a  natural elementwise multiplication $S_1 S_2 = \{s_1s_2: s_i \in S_i\}$,  making $ \mathcal P (\mathcal H  )$ a pre-semiring pair over $\tTz = \mcH$ under the action $aS = \{as :s\in S\}$.\footnote{Thus  $\zero a =\zero $ for all $a\in\mathcal H .$ Although not necessarily distributive, $\mathcal P (\mathcal H  )$       satisfies $(\boxplus _i S_i)(\boxplus _j S_j')\subseteq \boxplus _i (S_iS_j')$, cf.~\cite[Proposition 1.1]{Mas}.}
\item   A  hypergroup $\mathcal H$ is a \textbf{hyperfield} if $\mathcal H \setminus \{\zero\}$ is a  multiplicative group.
  \end{enumerate}
 \end{defn}

 \begin{rem}[\cite{AGR1,Row21}] \label{hp1}  Although \cite{NakR} uses hypermagmas, we work mostly with hypersemigroups and write $\mathcal{H}$ for  $(\mathcal{H},\boxplus,\zero)$.
 \begin{enumerate}\eroman
     \item  If it exists, the hyperzero 
     of a hypersemigroup is unique, by the familiar elementary argument.


   \item  $\mathcal{H}$ gives rise to a $\mcH$-module  $\mathcal{P}(\mcH),$ with   hyperaddition  given by $$S_1 \boxplus S_2 = \cup \{s_1 \boxplus s_2: s_i\in S_i\}.$$

 \item $\mathcal{P}(\mcH)$ is a weakly admissible $\tT$-bimodule in case
  \begin{enumerate}
      \item
  $\tT = \{\one\}$.

  \item    $\tT$ is a submonoid of $\mcH$, in particular when $  \mcH \setminus \{\zero\}$ is a monoid.
  \end{enumerate}

    \item Take any $\tT$-submodule $S_0$ of $\mcA:=\mathcal{P}(\mcH)$.  Then we get a $\tT$-pair $(\mathcal{P}(\mcH),\mathcal{P}(\mcH)_0)$ where
   $\mathcal{P}(\mcH)_0 = \{S \subseteq \mcH : S_0 \cap S \ne \emptyset\}. $

    \item Take  $S_0 = \{\zero\}$ in (iv). Thus  $\mathcal{P}(\mcH)_0 = \{S \subseteq \mcH : \zero \in S\},  $ the standard definition.
  \end{enumerate}\end{rem}

   \begin{definition}\label{hypp}    In~Remark~\ref{hp1}(v),
       the sub-pair of the pair
$(\mathcal{P}(\mcH),\mathcal{P}(\mcH)_0)$ generated by $\mcH$ is called the  \textbf{hyperpair} of $\mcH.$
    \footnote{This is the definition in \cite{AGR2}.  It could  lead to complications if $\mathcal{P}(\mcH)$  is not distributive, because then the  submagma of the power set $\mathcal{P}(\mcH)$ spanned by $\mcH$    need not be closed under multiplication.
}
   \end{definition}

\begin{rem}\label{mora}$ $
    \begin{enumerate}\eroman
  \item    Any hyperpair   has the  important surpassing relation $\subseteq,$
    i.e., $S_1 \preceq S_2$ when $S_1\subseteq S_2.$ Of course $\mathcal{P}(\mcH)$ is a lower semilattice with respect to $(\subseteq)$,  but not $\nsets(\mcH).$
      \item   More generally, if $\mcH$ has a  surpassing relation $\preceq$, then $\mathcal{P}(\mcH)$ has the surpassing relation given by  $S_1\preceq S_2$ if
      for each $s_1\in S_1$ there is $s_2\in S_2$ for which $s_1 \preceq s_2.$

    \item For hypersemigroups $\mcH$, the condition that $b_1\subseteq b_2$ implies $f(b_1) \subseteq  f(b_2)$ is automatic on $\mathcal{P}(\mcH)$.

   \item If each $a\in \mcH$ has a unique hypernegative, then $(\mcA,\mcA_0)$ has a negation map given by applying the hypernegative element-wise.
 \end{enumerate}
\end{rem}

   \begin{lem} Say a subset $S \subseteq \mathcal{P}(\mcH)$ is \textbf{weakly neutral} if $a\in a\boxplus S$ for each $a\in\mcH$. The  weakly neutral subsets comprise an $\mcH$-submagma of $\mathcal{P}(\mcH)$.
\end{lem}
\begin{proof} If $S_i$ are weakly neutral elements, then for any $a\in \mcH,$
    $a \boxplus (S_1\boxplus S_2) = (a\boxplus S_1)\boxplus S_2$, which contains $a\boxplus S_2,$  which in turn contains $a.$
\end{proof}

\subsection{Examples of hypersemigroups and their pairs}

 \begin{example}\label{Examplesofq}
    Let us first recall some of the celebrated hyperfields (and their accompanying hyperpairs, cf.~Definition~\ref{hypp}), from \cite{V1,V2}, as elaborated in \cite{BB,AGR1}. In every case we take $\tT=\mcH,$ and the hypersum of two distinct nonzero elements is not a singleton.
    \begin{enumerate}\eroman

 \item The \textbf{tropical hyperfield}
  consists of the set $\mcH=\R \cup \{-\infty\}$, with $-\infty$ here
  as the zero element and $0$ as the unit, equipped
  with the addition $a\boxplus b = \{a\}$ if $a>b$,
  $a\boxplus b = \{b\}$ if $a<b$,
  and $a\boxplus a= [-\infty,a]$.

        \item  The \textbf{hyperfield of signs} $L := \{0, 1, -1 \}$  has the intuitive multiplication
law, and hyperaddition defined by $1 \boxplus  1 = 1 ,$\ $-1
\boxplus  -1 = -1  , $ $\ x \boxplus  0 = 0 \boxplus  x = x $ for
all $x ,$ and $1 \boxplus  -1 = -1 \boxplus  1 = \{ 0, 1,-1\} .$

       \item  The \textbf{phase hyperfield}. Let $S^1$ denote the complex
  unit circle, and take $\mcH = S^1\cup \{ 0 \}$.
  Nonzero points $a$ and $b$ are \textbf{antipodes} if $a = -b.$
Multiplication is defined as usual (so corresponds on $S^1$ to
addition of angles). In this example, we denote  an open arc of less than 180 degrees
connecting two distinct points $a,b$ of the unit circle by
$({a, b})$. The hypersum is given, for $a,b\neq 0$, by
$$a \boxplus b=
\begin{cases} ({a, b}) \text{ if } a \ne \pm b  ;\\
  \{ -a,0,a \} \text{
if } a = -b \enspace , \\   \{  a \} \text{ if }   a = b\enspace .
\end{cases}$$
$\langle \tT \rangle$ contains only the points of $\mcH$, the elements of
the form $a\boxplus b$ with $a,b\in S^1$,
 and subsets  $C$ of $\mcH$,
where either $C$ is an open half circle connecting an element $a\in
S^1$ to $-a$ (which is obtained as the sum $a\boxplus b\boxplus
(-a)$, where $b \in C$), or $C=\mcH$ (which is obtained as $a\boxplus b\boxplus (-a)\boxplus (-b)$).

$e = \one \boxplus -\one = \{ -\one,0,\one \} ,$ so $e\boxplus e=e $. Thus
the {phase hyperfield} is    $\mcA_0$-idempotent, but is not $\circ$-distributive since   $(a,b)e = (a,b)\cup \{\zero\} \cup (-a,-b),$ whereas $(a,b)\boxplus (-a,-b)$ is the whole circle.

\item (The \textbf{weak phase hyperfield} \cite{AGR1};
  Akian-Gaubert's modification of the phase hyperfield).
  The non-zero elements  can still
  be represented by elements of the unit circle $S^1$. However,
  the hyperaddition differs -- for $a,b\neq 0$, we now have
  $$a \boxplus b=
  \begin{cases} {[a, b]}
    &
    \text{ if } a \ne b \text{ and } a\neq -b ,\\
  S^1 \cup \{0\} & \text{
    if } a = -b
  , \\   \{  a \} & \text{ if }   a = b .
\end{cases}$$
where ${[a, b]}$ denotes the {\em closed} arc of
less than $180$ degrees joining $a$ and $b$ (compare with the open
arc $(a, b)$ in the phase hyperfield).
Like the phase hyperfield, the weak phase hyperfield is
not distributive.

\item One can further modify the weak phase hyperfield, to get a  hypersemigroup which is not a hypergroup by defining $a\boxplus a= \emptyset.$

    \end{enumerate}
 \end{example}

\begin{example} Let $\mcH$ be a  set with an element stipulated $\zero.$ In each of the following cases, $\mcH$ will be a hypersemigroup whose hyperpair $(\mcA,\mcA_0)$   will satisfy property N.
   \begin{enumerate}\eroman  \item The following are uniquely negated, with $\mcA_0  $ as in Remark~\ref{hp1}(v).
     \begin{enumerate}
           \item  Define $\boxplus$ on $\mcH$ by $s \boxplus s = \{ \zero, s\}$, $\zero\boxplus s =s\boxplus \zero =s,$ and $s\boxplus s' =\mcH$ for all nonzero $s\ne s'$ in~$\mcH.$  Here $\mcA = \mcH  \cup \{\mcH\},$ and $-s=s.$
         \item  (As in \cite{Mas}) When $|\mcH|\ge 3,$ define $\boxplus$ on $\mcH$ by $s \boxplus s = \mcH\setminus s$ and $s\boxplus s' =\{ s,s'\}$ for all nonzero $s\ne s'$ in $\mcH.$    Here $\mcA = \mcH \cup_{s,s'\in \mcH} \{s,s'\} \cup_{s\in \mcH} (H\setminus \{s\}) \cup \{\mcH\}.$ Again $-s=s.$
      \item  When $|\mcH|\ge 4,$ define $\boxplus$ on $\mcH$ by $s \boxplus s' = \mcH\setminus \{s,s'\}$   for all nonzero $s, s' \in \mcH.$    Here  $\mcA = \mcH  \cup_{s,s'\in \mcH} \{s,s'\} \cup_{s\in \mcH} (H\setminus \{s\})\cup_{s,s'\in \mcH} (H\setminus \{s,s'\})  \cup \{\mcH\}.$
      (Here $H\setminus \{s\} = (s'\boxplus s'')\boxplus s)$ where $s',s''$ are distinct from $\{\zero,s\}.$)
     \end{enumerate}

          \item  Define $\boxplus$ on $\mcH$ by $s \boxplus s = s$ and $s\boxplus s' =\mcH$ for all $s\ne s'$ in $\mcH.$ This is   idempotent, and is not uniquely negated.

 \item    The following two examples are similar to those of \cite[Proposition~3.1]{NakR}. For all $a_i\in \mcH,$
\begin{enumerate}\eroman
      \item $a_1 \boxplus a_2 = \mcH$. This gives the weakly admissible pair  $ (\mcH   \cup \{\mcH\},\{\zero,\{\mcH\}\})$, which satisfies Property~N.
      \item $a_1 \boxplus a_2 = \emptyset$. $(\mcH\cup \{\emptyset\},\emptyset) $ is the minimal pair containing $\mcH$, but  is not a hyperpair.

 \end{enumerate}
     \item  $a_1 \boxplus a_2 = \{a_1, a_2\}$. Hence $S_1 \boxplus S_2 = S_1 \cup S_2$. $(\mathcal H,\mathcal H_0)$ fails Property N.

       \item Now take $\mcH$   ordered, with $\zero\in \mcH$ minimal.   Define $\boxplus$ on $\mcH$ by $s \boxplus s = \mcH$ and $s\boxplus s' =s'$ for all $\zero\ne s< s'$ in $\mcH.$ Taking $\mcA = \mcH   \cup \{\mcH\},$ and $\mcA_0 = \{\mcH\}$, the hyperpair $(\mcA,\mcA_0)$  is weakly admissible and $\mcA_0$-bipotent.


    \end{enumerate}
  In each case above, $\mcH$ becomes a hyper-semiring when $\mcH$ is a cancellative multiplicative monoid.
\end{example}

\subsection{Residue (quotient) hypermodules and hyperpairs}$ $

The following definition was inspired by Krasner~\cite{krasner}, a bit more general than \cite[Definition~3.1]{Row24}, but with the same verification as \cite[Lemma~3.2]{Row24}.
\begin{definition}\label{kras1} Suppose that  $\mcM$ is    a    $C$-bimodule where the ring $C$ is viewed as a multiplicative monoid $\tT$ as in Remark~\ref{cl0}(i), and $G$ is a normal multiplicative subgroup of $\tT $, in the sense that $bG = Gb$ for all $b\in \mcM$. Define the \textbf{residue hypermodule} $\mcH =\mcM/G$ over $\tT/G$
    to have multiplication induced by the cosets, and \textbf{hyperaddition} $\boxplus : \mcH \times \mcH \to \mathcal{P}(\mcH) $ by $$b_1 G \boxplus b_2 G = \{ cG: c\in b_1G + b_2G\}.$$

    When $\mcM = C$ is a field,   the residue  hypermodule $\mcH$ is a hyperfield (under multiplication of cosets), called the \textbf{quotient hyperfield} in the literature.
\end{definition}

As defined in Remark~\ref{hp1}(iv),
$(\mcM,\mcM_0)$ is a  pair, by \cite[Proposition 7.21]{AGR1}. For the applications in the hyperfield literature one would take $\mcM$ to be a  field  with a multiplicative subgroup  $\tT$.

\begin{remark}\label{hp2} In the residue hypermodule $\mcH =\mcM/G$,
\begin{enumerate}\eroman
  \item $\one_\mcH = G.$
    \item  $e = \one_\mcH \boxplus (-)\one_\mcH = \{ g_1-g_2: g_i \in G\}.$
    \item According to Definition~\ref{Hyp00}(vii), $Se = \{ b_1 g_1 - b_2 g_2: b_i \in S, g_i \in G\}$
    \item In particular, $ee = \{ (g_1-g_2)g_3 + (g_4-g_5)g_6 :  g_i \in G\}.$
     \item $e\boxplus e = \{ (g_1-g_2)  -(g_3-g_4) :  g_i \in G\}.$
\end{enumerate}
\end{remark}

$\mcM/G$ need not be $\circ$-distributive, since the phase hyperfield is a counterexample, cf.~Examples~\ref{Examplesofq}(iii), but $\mcM/G$ does satisfy $ee = e\boxplus e$ by \cite[Lemma~3.5]{Row24}.

\begin{lem}
    Any  surpassing relation $\preceq$ on $\mcM$ induces a surpassing relation on $\mcM/G$.
\end{lem}
\begin{proof}
    As in Remark~\ref{mora}, define $b_1 G \preceq b_2 G $ when for each $g\in G$ there is $g'\in G$ such that $b_1g \preceq b_2 g'$. We claim that if $b_iG \preceq b_i' G $ for $i = 1,2,$ then $b_1G \boxplus b_2 G \preceq b_1'G \boxplus b_2'G.$

    Indeed, if $b = b_1 g_1 + b_2 g_2 \in b_1G \boxplus b_2 G , $ then taking $b_i'g_i' \succeq b_ig_i$ in $b_i' G,$ we have $$b  = b_1 g_1 + b_2 g_2 \preceq   b_1'g_1' + b_2 'g_2' \in  b_1'G \boxplus b_2'G.$$
\end{proof}
\begin{example}
    \label{Examplesofquotient} A huge assortment of examples of quotient hyperfields is given in \cite[\S 2]{MasM}\footnote{In \cite{Row24} more examples are obtained when addition and  multiplication are reversed.}. Here are some of them.
We shall take $\mcH = \mcM/G,$ and its hyperpair $(\mcA,\mcA_0),$ as in Remark~\ref{hp1}(iv).  \begin{enumerate}\eroman
\item
$G = \{ \pm 1\}$. Then $\zero \in \mcH \boxplus \mcH ,$ so $(\mcA,\mcA_0)$ has characteristic 2 and is multiplicatively idempotent.
  \item  The {Krasner hyperfield} is  $F/F^\times$, for any field $F$.

\item  The {sign hyperfield} is $\R/\R^+.$

\item  The phase hyperfield can be identified  with the quotient hyperfield
$\C/\R_{>0}$.

\item  The weak phase hyperfield can be obtained by taking
  the quotient $\mcF/G$, where $\mcF=\C\{\{t^{\R}\}\}$, and $G$
  is the group of (generalized) Puiseux series with  positive real
  leading coefficient, where the leading coefficient is the coefficient
  $f_\lambda$ of the series $f=\sum_{\lambda \in \Lambda} f_\lambda t^\lambda$
  such that $\lambda$ is the minimal element of
  $\{\lambda \in \Lambda : f_\lambda \neq 0\}$.
\end{enumerate}

\end{example}

  \section{Tensor products of bimodules}

In \cite{JuR1}  tensor products were treated from the point of view of semirings, following \cite{Ka1,Tak}. There is a delicate issue here, concerning which are the morphisms in our category.
One may start with the more straightforward instance of homomorphisms. Then it is rather easy to construct a tensor product satisfying the theory of \cite{Jag,Po}, whose theorems are applicable.

When we bring in weak morphisms on both sides, things become much more complicated, and our approach   becomes more compatible with \cite{CHWW} and  \cite{NakR}.

 \subsection{Construction of the tensor product}
  \begin{example} The \textbf{free abelian semigroup} $(\mcS(X),+)$ on a set $X$ is the set of formal sums without parentheses, defined  by
         $X\subseteq \mcS(X)$ of length 1, and inductively, for $v,w \in \mcS(X)$ of lengths $m$ and $n$,  $v+w = w+v\in \mcS(X)$ of length $m+n.$
 \end{example}

 We follow the exposition of the classical tensor product, as in~\cite[Chapter~18]{Row08}. Throughout this section,   $\mathcal M_1$
 is a right $\tT $-module
 and $\mathcal M_2$ is a   left  $\tT $-module, for convenience\footnote{More generally, for $\mcM_1$  a right $\tT$-magma and $\mcM_2$   a left $\tT$-magma, we can define the \textbf{free magma} $(\mcF(X),*)$ on a set $X$ to be the set of formal expressions with parentheses, defined  by
         $X\subseteq \mcF(X)$ of length 1, and inductively, for $v,w \in \mcF(X)$ of lengths $m$ and $n$,  $(v *w)\in \mcF(X)$ of length $m+n,$ and the \textbf{$\tT$-tensor product magma} $\mathcal M_1  \otimes _{\tT} \mathcal M_2$ to be the magma $\mathcal F(\mcM_1 \times \mcM_2)/\Cong,$
where $\Cong$ is the
congruence
 generated by   $\big(\big( v_1*w_1,x_2\big),  \big( v_1, x_{2})*(w_1,  x_{2}\big)\big) ,\
 \big(\big( x_{1},   v_2*w_2\big),    (x_{1},v_2 ) *(x_{1},w_2 )  \big) ,$
 $  \big((  x_1 a, x_2 ), (x_1,a x_2
)\big)$
  $ \forall x_{i}, v_i,w_i\in \mathcal \mcM_i, $ $   a \in
\tT$. The remainder of this subsection would go over, mutatis mutandis, but the notation is more involved.}.

\begin{definition} \label{tp1} Let $\mathcal S(\mcM_1 \times \mcM_2)$ denote the free semigroup with base $\{ (v,w) \in \mcM_1 \times \mcM_2\}.$
\begin{enumerate}\eroman
  \item Define the \textbf{$\tT$-tensor product semigroup} $\mathcal M_1  \otimes _{\tT} \mathcal M_2$ to be the semigroup $\mathcal S(\mcM_1 \times \mcM_2)/\Cong,$
where $\Cong$ is the
congruence
 generated by all \begin{equation}\label{defconga}\bigg(\big( v_1+w_1,x_2\big),  \big( v_1, x_{2})+(w_1,  x_{2}\big)\bigg) ,\quad \bigg(\big( x_{1},   v_2+w_2\big),   (x_{1},v_2 ) +(x_{1},w_2 )  \bigg) ,\end{equation}
\begin{equation}
    \label{defcong1a} \bigg((  x_1 a, x_2 ), (x_1,a
 x_2
),\bigg)
\end{equation}
  $ \forall x_{i}, v_i,w_i\in \mathcal \mcM_i, $ $   a \in
\tT$.

      \item When $\mathcal M_1$
 and $\mathcal M_2$ have negation maps, we   incorporate the negation map into the definition of tensor product, by enlarging $\Cong$ to include also
 $ (\big((-) v_1,v_2\big),  \big( v_1, (-)v_2)\big) $ for all $v_i\in \mcM_i.$

\item   A \textbf{simple tensor} of $\mcM_1\otimes \mcM_2$ is an element of the form $ v_1\otimes  v_2$ for $v_i\in \mcM_i.$

 \end{enumerate}
 \end{definition}

Note that  this matches the situation that $\tT$ is a semiring, in which case we would just forget its addition.

\begin{rem}
    \label{ce2}$ $\begin{enumerate}\eroman
        \item
 The  sum of two non-simple tensors could be simple, modulo associativity. For example, if $v_3 =v_1 + v_2 $  then $(v_3 \otimes v_2)+(v_2 \otimes v_1) $ and $(v_2 \otimes v_2
)+(v_1 \otimes v_3) $ are non-simple   whereas
\begin{equation}
    \begin{aligned}
    (v_3 \otimes v_2)+((v_2 \otimes v_1) +(v_2 \otimes v_2
))+(v_1 \otimes v_3)) &= (v_3 \otimes v_2)+((v_2 \otimes v_3
))+(v_1 \otimes v_3))\\& =  (v_3 \otimes v_2)+(v_3\otimes v_3)  =  v_3 \otimes (v_2+v_3) .\end{aligned}\end{equation} 
\item
     If $a_1 v_2 = a_2 w_2,$ then $v_1 a_1 \otimes v_2 + w_1 a_2\otimes w_2 = (v_1+w_1)\otimes a_2 w_2$.
    This process of ``recombining'' is straightforward when $\tT$ is a group, since then it is applicable  for any sum
    $$v_1 \otimes v_2 + w_1\otimes w_2 = (v_1a_1 ^{-1}+w_1  a_2^{-1})\otimes a_2 w_2. $$
    \end{enumerate}
 \end{rem}

 \begin{rem}\label{tpm} As in \cite{CHWW}, one could define the tensor product   $\tT_1 \otimes_\tT \tT_2$ of a right $\tT$-monoid  $\tT_1$ and  a left $\tT$-monoid
$\tT_2$ to be $(\tT_1\times \tT_2)/\Phi$ where  $\Cong$ is the
congruence
 generated by all $\bigg((  x_1 a, x_2 ), (x_1,a x_2
)\bigg)$,
  $ \forall x_{i}\in \mathcal \mcM_i, $ $   a \in
\tT$.
 \end{rem}

\subsubsection{Tensor products of bimodule prepairs}$ $

In order to define the tensor product  bimodule prepair, we need to define $(\mcM_1\otimes \mcM_2)_0$.

\begin{definition}  Suppose $(\mcM_1,{\mcM_1}_0)$ is a $(\tT_1,\tT)$-bimodule prepair and $(\mcM_2,{\mcM_2}_0)$ is a $(\tT,\tT_2)$-bimodule  prepair.  Define
$(\mcM_1\otimes \mcM_2)_0 := ({\mcM_1}_0\otimes_{\tT} \mcM_2)+ (\mcM_1\otimes_{\tT} {\mcM_2}_0).$
\end{definition}

\begin{lem}
    $((\mcM_1\otimes \mcM_2) ,(\mcM_1\otimes \mcM_2)_0) $ is a $(\tT_1,\tT_{2})$-bimodule prepair.
\end{lem}
\begin{proof}
       $((\mcM_1\otimes \mcM_2) ,(\mcM_1\otimes \mcM_2)_0) $ is closed under left multiplication by $\tT_1$ and right multiplication by~$\tT_2.$ $(\mcM_1\otimes \mcM_2)_0$ is closed under addition.
\end{proof}

 \subsubsection{Pre-orders on the tensor product}$ $

 Note that the multiplicative structure in the tensor product could be carried out along the lines of \cite{CHWW}, but coping with addition makes the situation much more intricate.

 \begin{lem}
     If $\mcM_i$ each have a pre-order $\preceq$, then $\mcM_1 \otimes \mcM_2$ has the pre-order given by $x \preceq y$ if whenever we can write $x =\sum v_i\otimes w_i$ then we can write $y =\sum v'_i\otimes w'_i$ such that  $ v_i \preceq v'_i $ and $  w_i \preceq   w_i'$ for each $i.$
 \end{lem}
\begin{proof}
    Clearly this condition is  transitive and passes to sums.
\end{proof}

\begin{rem}
    In  general, for semiring pairs, one conceivably might have $\one \otimes \one \in (\mcM_1\otimes \mcM_2)_0.$ In such a case the tensor product of surpassing relations cannot be a surpassing relation. But one must realize that even in the classical situation we could have $\one \otimes \one = \zero.$
\end{rem}

 \subsection{Properties of tensor products}\label{tprop}$ $

\begin{definition}\label{bal1}$ $
    \begin{enumerate}\eroman
    \item A \textbf{1-balanced map} is a function $\psi: \mathcal M_1  \times   \mathcal M_2 \to \mcM$ where $\mcM$ is a semigroup, satisfying
    \begin{equation}
    \psi(v_1a, v_2) =\psi(v_1, av_2),
\end{equation} for all   $v_i \in \mcM_i,$ $a\in \tT.$

\item A \textbf{balanced map} is a 1-balanced map satisfying
\begin{equation}
    \psi(v_1 +w_1, x_2) =\psi(v_1 ,x_2)+\psi(w_1, x_2),
\end{equation}

  \begin{equation}
    \psi(x_1, v_2 +w_2 ) =\psi(x_1,v_2 )+\psi(x_1,w_2),
\end{equation}
for all   $v_i, w_i, x_i\in \mcM_i.$
    \end{enumerate}

    When $\mathcal M_1, \mathcal M_2$  are given with negation maps, then a balanced map is required to satisfy $((-)v_1)\otimes v_2 = v_1 \otimes (-v_2),$
    which we define to be $(-)(v_1\otimes v_2)$.
\end{definition}

\begin{thm}\label{tenf} Given $(\tT_1,\tT)$-bimodules $\mcM_1, \mcN_1$,   and  $(\tT,\tT_2)$-bimodules $\mcM_2, \mcN_2,$ suppose $f_1:\mcM_1\to \mcN_1$  is a   $(\tT_1,\tT)$-bimodule homomorphism              and   $f_2:\mcM_2\to \mcN_2$  is a  $(\tT,\tT_2)$-bimodule homomorphism.  Then there is a $(\tT_1,\tT_2)$-bimodule          homomorphism $$f_1\otimes f_2 :  \mathcal M_1  \otimes _{\tT} \mathcal M_2 \to   \mathcal N_1  \otimes _{\tT} \mathcal N_2 $$ given by $(f_1\otimes f_2)(v_1 \otimes v_2) =f_1(v_1)\otimes f_2(v_2),$ which sends $(\mathcal M_1  \otimes _{\tT} \mathcal M_2 )_0 $ to $( \mathcal N_1  \otimes _{\tT} \mathcal N_2)_0.$

\end{thm}

\begin{proof}
    The balanced map $\Psi : (v_1, v_2) \mapsto f_1(v_1)\otimes f_2(v_2)$  sends $\mathcal M_1  \otimes _{\tT} \mathcal M_2  $ to $   \mathcal N_1  \otimes _{\tT} \mathcal N_2$ seen by applying $f_1$ and $f_2$ to the respective components.

    To see that $\mathcal M_1  \otimes _{\tT} \mathcal M_2  $ is a left $\tT_1$-module, take $f_1$ to be left multiplication by $a$ and $f_2$ to be the identity map in the previous paragraph. Symmetrically, one sees that  that $\mathcal M_1  \otimes _{\tT} \mathcal M_2  $ is a right $\tT_2$-module, and it is seen to be a bimodule by taking $f_1$ to be   left multiplication by $a_1$, and $f_2$ to be right multiplication by $a_2$.

    $(f_1\otimes f_2)(\mcM_1\otimes \mcM_2)_0 = f_1({\mcM_1}_0)\otimes_{\tT} f_2(\mcM_2)+ f_1(\mcM_1)\otimes_{\tT} f_2({\mcM_2}_0) \subseteq {\mcN_1}_0\otimes_{\tT} \mcN_2 +  \mcN_1 \otimes_{\tT}  {\mcN_2}_0. $
\end{proof}

\begin{cor}\label{tensorisoa} $ $ Assume throughout that $\mcM_1$ is a   $(\tT_1,\tT)$-bimodule and  $\mathcal M_2$ is a  $(\tT,\tT_2)$-bimodule.
      \begin{enumerate}\eroman

\item   If $\mcA_1,\mcA_2$ are weakly admissible, then   $\mathcal M_1
\otimes \mathcal M_2$ are weakly admissible over $\one\otimes \tT = \tT \otimes \one$.

 \item Any homomorphism $f:\mcM_1 \otimes \mcM_2 \to \mcN$ is determined by its action on simple tensors.

  \item If both $\mcM_1$ and $\mcM_2$ are  $\tT$-bimodules, then $\mcM_1 \otimes_\tT \mcM_2 \cong \mcM_2 \otimes_\tT \mcM_1 $ via $(v_1 \otimes v_2)\mapsto v_2\otimes v_1.$


    \end{enumerate}
\end{cor}
 \begin{proof} Repeated applications of Theorem~\ref{tenf}.

     (i) Identifying $\tT_i \subseteq \mcM_i,$ enables us to identify $\tT_1\otimes 1$ inside $\mcM_1 \otimes \mcM_2.$

    (ii) $f$ is determined by the balanced map  $\psi:\mcM_1 \times \mcM_2 \to \mcN$ given by $\psi(v_1,v_2) = f(v_1\otimes v_2).$

(iii) Take  the balanced map $(v_1,v_2)\mapsto v_2\otimes v_1$.

 \end{proof}

\begin{rem}\label{ce1}$ $ \begin{enumerate} \eroman
    \item  If one takes $\tT = \{\one\}$ as in  Remark~\ref{hp1}(iii), then Equation~\eqref{defcong1a} is redundant.
     \item Associativity of addition in $\mcM_1 \otimes \mcM_2$ follows from associativity of addition in $\mcS(\mcM_1 )\times \mcS(\mcM_2 ).$

\end{enumerate}
\end{rem}

\begin{cor}\label{tensoriso0}
    The tensor product in this subsection satisfies the criteria of \cite[Definitions 1,2]{Po}.
\end{cor}
\begin{proof}
    By Corollary~\ref{tensorisoa}.
\end{proof}

 \begin{cor}\label{tensoriso} Suppose $\mcM,\mcM_i$  are $\tT_1,\tT$-bimodules and
 $\mcN,\mcN_j$  are $\tT,\tT_2$-bimodules.
    \begin{enumerate}\eroman   \item $\tT_1 \otimes _{\tT_1} \mcM \cong \mcM$ and  $  \mcM \otimes _{\tT} {\tT} \cong \mcM$,
        \item   $(\oplus \mcM_i) \otimes N \cong \oplus (\mcM_i \otimes N)$,
  \item  $ \mcM  \otimes (\oplus \mcN_i) \cong \oplus (M \otimes \mcN_i)$,

   \item  $ (\oplus \mcM_i)  \otimes (\oplus \mcN_j) \cong \oplus (\mcM_i \otimes \mcN_j)$,
   as $\tT_1,\tT_2$-bimodules.

    \end{enumerate}
 \end{cor}
 \begin{proof} These are results in \cite{Jag,Po}.
 \end{proof}

 \begin{cor}\label{tensoriso1} Suppose $\mcM_i$  are $\tT_i,\tT_{i+1}$-bimodules. Then $(\mcM_1 \otimes \mcM_2)\otimes \mcM_3 \cong \mcM_1 \otimes (\mcM_2\otimes \mcM_3) $ as $\tT_1,\tT_3$-bimodules.
 \end{cor}
\begin{proof} The proof is   from \cite[Theorem~2]{Po}.
\end{proof}

\begin{prop}
    If $\mcM_i = \mathcal P(\mcH_i)$ where   $\mcM_1$ is a right $\tT$-module and  $\mathcal M_2$ is a left $\tT$-module, then $\mcM_1 \otimes \mcM_2$ is $\subseteq$-distributive.
\end{prop}
\begin{proof}
    We shall show that $S\otimes (S_1'\boxplus S_2')\subseteq S S_1'  \, \boxplus \, SS_2' $ for all $S\subseteq \mcH_1$ and $S_1',S_2' \subseteq \mcH_2$.  By definition, $S\otimes (S_1'\boxplus S_2') = \{ a \otimes  a' : a\in S, \ a' \in S_1' \boxplus S_2'\} = \{a (S_1'\boxplus S_2'): a\in S\}.$
    Thus we need to show that $a \otimes (S_1'\boxplus S_2') \subseteq  a  \otimes S_1'\, \boxplus \, a  \otimes S_2'$ for each $a\in S.$ But for each $a'\in S_1'\boxplus S_2'$, if
    $a' \in b_1'\boxplus b_2'$ for suitable nonzero $b_i'\in S_i'$, then  $a\otimes a' \in ( a \otimes  b_1') \boxplus (a\otimes b_2'), $ and we are done by induction on height. Hence we may assume that $a' \in \mcH_2,$ and the assertion is immediate.
\end{proof}

Write $\tT^\op$ for the opposite monoid. Then
any $\tT_1,\tT_2$ bimodule  is a $\tT_1\times \tT_2^\op$ module, under the operation $(a_1,a_2)b = a_1 b a_2,$ by the obvious verification.
\begin{lem}\label{tp4}
    Suppose $\tT_i \supseteq \tT$ are commutative monoids. Then   $\mcM_1\otimes \mcM_2 $ is a $\tT_1\otimes_\tT \tT_2$-module,  under the action $(a_1\otimes a_2) b = a_1 b a_2.$
\end{lem}
\begin{proof}
    We have the $\tT_1\times \tT_2$ action, which induces a $\tT_1\otimes_\tT \tT_2$-action since $(a_1a\otimes a_2) b = a_1a b a_2 = a_1 b aa_2 = (a_1\otimes aa_2) b .$
\end{proof}

\begin{cor}
     If $\mcA_i$ are $\tT_i$-semirings with each $\tT_i\subseteq \tT$ and if $\mcA_1 \otimes \mcA_2$ is a semiring as in Corollary~\ref{tensoriso}(vi), then $\mcA_1\otimes\mcA_2$ is a $\tT_1 \otimes _\tT \tT_2$-semiring.
\end{cor}

\subsection{Residue tensors viewed functorially}$ $

One nice feature of this approach is that the residue functor preserves tensor products.
The notation for hyperaddition in residue hypermodules continues to be $\boxplus.$

\begin{prop}\label{res1}
    Suppose that $\tT _1, \tT_2$ are abelian monoids, $\mcM_i$ are      $\tT_i$-modules for $i =1,2,$ and $G$ is a common     subgroup of $\tT _1$ and $\tT_2.$
    The   tensor product  of the   residue hypermodules $\mcH_1 =\mcM_1/G_1$ and $\mcM_2 =\mcM_2/G_2$ is isomorphic to $(\mcM_1 \otimes \mcM_2)/ (G_1 \otimes G_2)$.
\end{prop}
\begin{proof}
 The natural   map  $\mcM_1 \otimes \mcM_2 \to ( \mcM_1/G_1 )\otimes (\mcM_2/G_2) $ given by $y_1 \otimes y_2 \mapsto y_1G_1 \otimes y_2G_2  $ preserves the defining congruence, thereby inducing a map of the cosets, which is 1:1 and onto.
\end{proof}

 \subsection{Tensor   extensions}\label{tensex1}$ $

The theory works most smoothly when  the left side of the tensor product is a monoid, since we are concerned only with addition on the right.

\begin{definition}\label{tensext} Suppose $\mcM$ is a $\tT$-module, and $\tT'$ is a monoid containing $\tT$.
\begin{enumerate} \item
   Define the \textbf{$\tT$-tensor extension}  $ \tT'  \otimes _{\tT} \mathcal M$ to be the semigroup $\mcS(  \tT' \times   \mcM)/\Cong,$
where $\Cong$ is the
congruence
 generated by all \begin{equation}\label{defcongaa}  \bigg(\big( a',   v+w\big), \big( a',   v \big)+\big( a',    w\big)\bigg)\qquad   \bigg((  a' a, w ), (a',a
 w
)\bigg) \qquad \forall   a,  a' \in
\tT,\ v ,w \in \mathcal \mcM.
\end{equation}

     When
   $\mathcal M$ has a negation map, we   incorporate the negation map into the definition of tensor product, by enlarging $\Cong$ to include also
 $ (\big((-) a,w),  \big( a, (-)w   )\big) $ for $a\in \tT,$ $w\in \mcM.$


 \end{enumerate}
 \end{definition}

\begin{rem}
    One sees easily that if $\mathcal M$ is   weakly admissible, then so is   $ \tT'  \otimes _{\tT} \mathcal M$.
\end{rem}

Now we modify \S\ref{tprop}, starting with Definition~\ref{bal1}.

\begin{definition}\label{bal2}$ $
    \begin{enumerate}\eroman
    \item A \textbf{balanced map} is a function $\Psi: \tT'  \times   \mathcal M \to \mcN$ where $\mcN$ is a semigroup, satisfying
\begin{equation}
    \Psi(a' , w_1 +w_2 ) =\Psi(a' ,w_1 )+\Psi(a ',w_2),
\end{equation}
\begin{equation}
    \Psi(a'a, w_1) =\Psi(a', aw_1),
\end{equation} for all  $a\in \tT,$  $a'\in \tT',$  $v_i, w_i \in \mcM_i.$
    \end{enumerate}
\end{definition}
 \begin{lem}
 $\mcM' := \tT'  \otimes _{\tT} \mathcal M$ is a $\tT'$-module.
 \end{lem}
 \begin{proof}
   The map $\tT' \times \mcM' \to \mcM' $ given by $(a,\sum a_i\otimes w_i)\mapsto\sum a a_i\otimes w_i$ factors through the congruence~$\Cong.$
 \end{proof}

As in Theorem~\ref{tenf}, we have the key observation leading to the isomorphism theorems of Corollaries`S\ref{is1} and  \ref{is2} below.

\begin{thm}
    \label{tenfa} $ $\begin{enumerate}\eroman
\item Suppose $\mcM_i$ are $\tT$-modules for $i=1,2$.
If $f_1:\tT_1\to \tT_2$  is a   monoid homomorphism   fixing $\tT$           and   $f_2:\mcM_1\to \mcM_2$  is a  module homomorphism,  then there is a $\tT$-module          homomorphism $$f_1\otimes f_2 :  \tT_1  \otimes _{\tT} \mathcal M_1 \to  \tT_2  \otimes _{\tT} \mathcal M_2 $$ given by $(f_1\otimes f_2)(a \otimes w) =f_1(a)\otimes f_2(w).$

 \item If $\mcM$ is a $\tT$-semiring and $\tT_1$ is commutative, then $\tT_1 \otimes _\tT \mcM$ is  a $\tT_1$-semiring.

 \item Suppose $\tT_1, \tT_2$ are monoids containing $\tT$, and $\mcA_1, \mcA_2$ are $\tT$-semirings.
If $f_1:\tT_1\to \tT_2$  is a   monoid homomorphism              and   $f_2:\mcA_1\to \mcA_2$  is a semiring homomorphism,  then there is a $\tT_1$-semiring        homomorphism $$f_1\otimes f_2 :  \tT_1  \otimes _{\tT} \mathcal A_1 \to  \tT_2 \otimes _{\tT} \mathcal A_2 $$ given by $(f_1\otimes f_2)(a \otimes b) =f_1(a)\otimes f_2(b),$ for $a\in \tT_1,\ b\in \mcA_1$.
\end{enumerate}
\end{thm}
\begin{proof}
 (i)   The balanced map $\Psi : (a, w) \mapsto f_1(a)\otimes f_2(w)$  sends the congruence $\Phi_1$ of Definition~\ref{tp1} to $\Phi_2,$ seen by applying $f_1$ and $f_2$ to the respective components.

 (ii) For $a_i \in \tT$ and $w_i\in \mcA_i,$ consider the balanced map $f_{a_2;w_2}\tT_1 \times \mcM \to \tT_1 \otimes \mcM$ given by $(a_1 ,w_1 ) \mapsto a_1 a_2 \otimes w_1w_2$.
 This induces a right multiplication map $\tT_1 \otimes \mcM_1 \to \tT_1 \otimes \mcM_1$.

 (iii) $f_1\otimes f_2$ respects multiplication, namely
 $$(f_1\otimes f_2)((a_1\otimes w_1)(a_2\otimes w_2)) = ((f_1\otimes f_2)(a_1\otimes w_1))((f_1\otimes f_2)(a_2\otimes w_2)).  $$
\end{proof}

The most important case is when $\tT_1 =\tT$ and $f_1:\tT \to \tT_2$ is the inclusion map.

\begin{cor} \label{is1}
 Let    $ \tT_1$ be a monoid containing $\tT$, and $\mcM_1, \mcM_2$ be $\tT$-modules.
If   $f:\mcM_1\to \mcM_2$  is a  bimodule homomorphism,  then there is a bimodule          homomorphism $$1\otimes f :  \tT_1  \otimes _{\tT} \mathcal M_1 \to  \tT_1  \otimes _{\tT} \mathcal M_2 $$ given by $(1\otimes f)(a \otimes w) =a\otimes f(w).$
\end{cor}

\begin{cor} \label{is2} Suppose that $\mcA$ is a $\tT$-semiring.
  \begin{enumerate}
     \eroman
      \item If $\tT \subseteq \tT_1 \subseteq \tT_2,$ then $\tT_2 \otimes (\tT_1 \otimes \mcA) \cong  \tT_2\otimes \mcA.$
         \item $\tT \otimes \mcA \cong \mcA$  as semirings.
            \item If $\tT \subseteq \tT_1$ and $\mcA_1 \to \mcA_2$ is an epimorphism of $\tT$-semirings, then  $\tT_1\otimes\mcA_1 \to \tT_1\otimes\mcA_2$ is an epimorphism.
  \end{enumerate}
\end{cor}
 \begin{proof}
  As in \cite[Theorem~3]{Po}.
 \end{proof}

 We also get a pair.

 \begin{thm}\label{tenfb} Suppose $\tT_1$ is   the disjoint  union $ \underset{\cdot}{\cup}\,
      c_i \tT$, and $\mcM$ is a $\tT$-module.  \begin{enumerate}\eroman
          \item  Every element of $\tT_1 \otimes \mcM$ can be written  in the form $\sum_{i\in I} c_i \otimes w_i$, for uniquely determined $w_i\in \mcM$.

          \item Define ${\tT_1 \otimes \mcM}_0 = \tT_1 \otimes \mcM_0.$ Then ${\tT_1 \otimes \mcM}_0 = \{ \sum c_i \otimes w_i : w_i\in \mcM_0\},$ and $({\tT_1 \otimes \mcM},{\tT_1 \otimes \mcM}_0)$ is a pair.
      \end{enumerate}
 \end{thm}
\begin{proof}
    (i)    The projection $\pi_j:\tT_1\to \tT$ onto the  coefficient of $c_j$   yields a  projection $1\otimes \pi_j:\tT_1 \otimes\mcM_2\to  \mcM_2$. If
     $\sum_{i\in I} c_i \otimes w_i= \sum_{i\in I} c_i \otimes w'_i$ then applying $\pi_j \otimes 1$ yields $ w_j=    w'_j.$

     (ii) If $\sum c_i a_i \otimes w_i \in \tT_1 \otimes \mcM_0$ for $w_i \in \mcM_0,$ then $\sum c_i a_i \otimes w_i =  \sum c_i \otimes  a_iw_i \in \tT_1 \otimes \mcM_0.$ The second assertion is clear.
\end{proof}

\begin{thm}\label{m2}
    Suppose $\tT_1, \tT_2$ are monoids containing $\tT$ where $\tT_1$ is the disjoint union $ \underset{\cdot}{\cup}\,
      c_i \tT$, and $\mcM_1, \mcM_2$ are $\tT$-modules.
If    $f:\mcM_1\to \mcM_2$  is a weak    morphism (resp.~$\preceq$-morphism),  then there is a weak  morphism (resp.~$\preceq$-morphism)$$ \tilde f:= 1\otimes f :  \tT_1  \otimes _{\tT} \mathcal \mcM_1 \to  \tT_2 \otimes _{\tT} \mathcal \mcM_2 $$ given by $ \tilde f(\sum c_i  \otimes w_i) =\sum c_i  \otimes f(w_i),$ for $ w_i\in \mcM_1$.
\end{thm}
\begin{proof}
    Applying Theorem~\ref{tenfb} gives us a projection $\pi_j \otimes 1 :  \tT_1  \otimes _{\tT} \mathcal \mcM_1 \to  \mcM_1$  given by $\sum _i c_ia_i \otimes b \mapsto a_j \otimes b = 1\otimes a_j b,$ and $1 \otimes f$ is well-defined.

    If $\sum_k \sum_i c_i \otimes w_{i,k} \in \{\tT_1 \otimes \mcM\}_0 $ for $w_{i,k} \in \mcM,$ then $ \sum_i \sum_k   c_i \otimes w_{i,k}= \sum _i (c_i \otimes \sum_k w_{i,k}),$ implying $\sum_k w_{i,k} \in \mcM_0,$ implying
\begin{equation}
    \begin{aligned}     \qquad
       \tilde f & (\sum_{i,k} c_i  \otimes w_{i,k}) =\sum_i\left( c_i  \otimes \sum_k f(w_{i,k})\right)   \in (\tT_2 \otimes _{\tT} \mathcal \mcM_2)_0,
    \end{aligned}
\end{equation}
proving the assertion for weak morphisms. The proof for $\preceq$-morphisms is analogous.
\end{proof}

 \section{Tensor products of   morphisms  which need not be homomorphisms} $ $

Now we get to the sticky part.
The proof of Theorem~\ref{tenf} fails for   morphisms  which are not homomorphisms, and thus its usefulness for hypersemigroups is limited. Our goal nevertheless is to salvage whatever information we can. We do not even have an obvious functor for residue hypersemigroups in general.
          \begin{rem}
Suppose  $f_i:\mcM_i \to \mcN_i$ are homomorphisms of $\tT_i$-modules, for $\tT_1$ and $\tT_2$ abelian, and $G$ is a common subgroup of $\tT_1$ and $\tT_2.$ Then writing $\overline{\mcM_i  }= \mcM_i/G$, $\overline{\mcN_i  }= \mcN_i/G,$ and  $\bar f_i: \overline{\mcM_i  } \to \overline{\mcN_i  }$ for the  weak morphism induced by $f_i$, we would want  a balanced map $\Psi:=\bar f_1\times \bar f_2:  {\overline{\mcM_1  }  }\times  {\overline{\mcM_2  }  }\to \overline{\mcN_1  }\otimes \overline{\mcN_2  }$  given by
    $$ (y_{1} \times y_{2}) \mapsto f_1(y_{1}) G f_2(y_{2}) = \{ f_1(y_{1}) g f_2(y_{2}) :g\in G\},$$
        to induce a
        weak morphism  $ \overline{\mcM_1  }\otimes \overline{\mcM_2  }\to \overline{\mcN_1  }\otimes \overline{\mcN_2  }$.

  But
  $       \Psi(v_1 +w_1, x_2) = f_1(v_1 +w_1)Gf_2(x_2)$ whereas $
         \Psi(v_1 ,x_2)\boxplus \Psi(w_1, x_2)= f_1(v_1)Gf_2(x_2) \boxplus f_1(w_1)Gf_2(x_2),$
which could be larger. For example, for $G = \{+1,-1\},$ the first term is  $\pm  f_1(v_1 +w_1)f_2(x_2) $ whereas the second term is $\pm(f_1(v_1)Gf_2(x_2) \pm f_1(w_1)Gf_2(x_2)).$ Thus, the tensor product of homomorphisms modulo the residue map is too naive in general.
     \end{rem}

\begin{rem}\label{NaR}
       In \cite{NakR}, to expedite a categorical approach for hypermagmas, the tensor product of $\mcM_1$ and~$\mcM_2$  is defined as satisfying
    \begin{equation}
        \label{defc1} (v_1\otimes v_2)\boxplus (v_1\otimes v_2) =(v_1\boxplus v_1)\otimes v_2 \ \cup\ v_1\otimes (v_2\boxplus v_2),
         \end{equation}
    \begin{equation}    (v_1\otimes v_2)\boxplus  (v_1'\otimes v_2) =(v_1\boxplus v_1')\otimes v_2,
    \end{equation}
     \begin{equation} \label{defc2}   (v_1\otimes v_2)\boxplus  (v_1\otimes v_2') = v_1\otimes (v_2\boxplus v_2'),
    \end{equation}
    \begin{equation} \label{defc3}   (v_1\otimes v_2)\boxplus  (v_1'\otimes v_2') = \emptyset,
    \end{equation}
for $v_1\ne v_1' \in \mcM_1$ and $v_2\ne v_2' \in \mcM_2$.
As they point out, this is the smallest hypermagma satisfying $ (v_1\boxplus v_1')\otimes v_2 \subseteq (v_1\otimes v_2)\boxplus  (v_1'\otimes v_2) $ and $ v_1\otimes (v_2\boxplus v_2') \subseteq (v_1\otimes v_2)\boxplus  (v_1\otimes v_2')  $
   for all   $v_1, v_1' \in \mcM_1$ and $v_2, v_2' \in \mcM_2$,
and has excellent functorial properties, but there is considerable collapsing, and  associativity of addition even fails   in the bipotent situation: \end{rem}

\begin{example}
    \label{NaR1}
If $v_1\boxplus v_2 = v_1$ and  $w_1\boxplus w_2 = w_1$ then
$((v_1 \otimes w_1) \boxplus (v_1 \otimes w_2)) \boxplus ((v_2 \otimes w_1) \boxplus (v_2 \otimes w_2))= (v_1 \otimes w_1)\boxplus (v_2 \otimes w_1)= (v_1 \otimes w_1) $
whereas
$$(v_1 \otimes w_1) \boxplus ((v_1 \otimes w_2) \boxplus (v_2 \otimes w_1)) \boxplus (v_2 \otimes w_2)) = (v_1 \otimes w_1) \boxplus (\emptyset \boxplus (v_2 \otimes w_2)) =\emptyset.$$
\end{example}

Thus, although \cite{NakR} is   appropriate   for $\preceq$-morphisms in the category of hypermagmas, we feel that the more traditional approach may be needed in the category of hypermodules.

\subsection{Defining $f_1\otimes f_2$ via simple tensors}$ $

Given morphisms $f_i : \mathcal M_i \to \mathcal N_i$ for $i = 1,2,$
we search for a well-defined morphism $f_1 \otimes f_2 : \mathcal M_1
\otimes \mathcal M_2 \to \mathcal N_1 \otimes \mathcal N_2$.

Unfortunately, there does not seem to be a single approach that works for everything.

\begin{example}$ $
    \begin{enumerate} \eroman
        \item One might be tempted to define $\tilde f = f_1\otimes f_2$ by $\tilde f(v\otimes w) = f_1(v) \otimes f_2(w)$ and $f_1\otimes f_2$ to be $\zero$ on non-simple tensors. However, this would lead to the situation where $\tilde f(\la_1 \otimes \la_2 +\la_2 \otimes \la_1) = \tilde f(\la_1 \otimes \la_1 +\la_2 \otimes \la_2) =\zero $ whereas $\tilde f((\la_1 + \la_2)\otimes (\la_1 + \la_2)) = f_1(\la_1 + \la_2)\otimes f_(\la_1 + \la_2).$

         \item Even worse, the decomposition of a tensor into the sum of two simple tensors need not be unique. For example
         $$\la_1  \otimes \la_2 +(\la_1 +\la_2) \otimes \la_3 = \la_1  \otimes(\la_2 +\la_3) + \la_2
 \otimes \la_3 .$$
          \end{enumerate}
\end{example}

Here is the piece that we can salvage in general.

\begin{prop}
    \label{tenf1} Suppose $\mcM_1, \mcN_1$ are   $\tT$-modules and $\mcM_2, \mcN_2$ are left $\tT$-modules.
If $f_1:\mcM_1\to \mcN_1$             and   $f_2:\mcM_2\to \mcN_2$ are   multiplicative  maps, then there is a 1-balanced map, which we denote as $f_1\times_\circ f_2 : \mcM_1 \times \mcM_2 \to \mcN_1 \otimes \mcN_2$,  given by   $(b_1,b_2)\mapsto f_1(b_1)\otimes f_2(b_2).$
\end{prop}

\begin{proof}
    $(f_1\times_\circ f_2) (b_1a,b_2) = f_1(b_1 )a \otimes f_2(b_2) =  f_1(b_1 )\otimes  a f_2(b_2) = (f_1\times_\circ f_2) (b_1,ab_2)$ for $a\in \tT.$
\end{proof}

 \begin{rem}\label{2act}
     If $\mcM_1, \mcN_1$ are   $(\tT_1,\tT)$-bimodules and $\mcM_2, \mcN_2$ are   $(\tT,\tT_2)$-bimodules, then $\mcM_1 \times \mcM_2$ is a $\tT_1,\tT_2$-bimodule, where the $\tT_1$ action is on the left on $\mcM_1,$ and the $\tT_2$ action is on the left on $\mcM_2.$
     Then $f_1\times_\circ f_2 $ of Proposition~\ref{tenf1}  preserves these left and right actions.
 \end{rem}
\begin{definition}

 Given a  $(\tT_1,\tT)$-bimodule prepair $( {\mcM_1}, {\mcM_1}_0)$,   a $\tT,\tT_2$-bimodule prepair   $( {\mcM_2}, {\mcM_2}_0)$, and   a $\tT_1,\tT_2$-bimodule prepair   $(\mcN,\mcN_0)$, define $\WMor (\mcM_1 \times \mcM_2,\mcN)$ to be the set of $\tT_1,\tT_2$ multiplicative maps $f:\mcM_1 \times \mcM_2 \to \mcN$ (as in Remark~\ref{2act}) satisfying  the property that if $b_j, b_{i,j} \in \mcM_i$ with  $\sum _j b_{i,j} \in {\mcM_i}_0,$ then $$f\left( \sum _j b_{1,j}, b_2\right)\, , \, f\left(b_1, \sum _j b_{2,j}\right)\in \mcN_0.$$
\end{definition}
\begin{lem}
      If $f_i:(\mcM_i,{\mcM_i}_0)\to (\mcN_i,{\mcN_i}_0)$ are weak morphisms,  then $$f_1 \times_\circ  f_2 \in \WMor (\mcM_1 \times \mcM_2,\mcN_1 \otimes \mcN_2).$$
\end{lem}
\begin{proof}
 $(f_1 \times_\circ  f_2)( \sum _j b_{1,j}, b_2) = f_1(\sum _j b_{1,j})\otimes f_2(b_2) \in {(\mcN_1 \otimes \mcN_2})_0.$ Likewise for the other side.
\end{proof}

Nevertheless, note that $\WMor (\mcM_1 \times \mcM_2,\mcN_1 \otimes \mcN_2)$ need not induce weak morphisms in the tensor product, since there may be some weird sum in $(\mcM_1 \times \mcM_2)_0$.

\subsubsection{Digression: Hypertensor products of maps}$ $

An intriguing (but problematic) approach   is to    introduce  a hyperstructure.

\begin{remark}\label{rev}
    Given maps $f_i : \mathcal M_i \to \mathcal N_i$, $i = 1,2,$    define
 $ f:=f_1 \otimes_\supseteq  f_2 : \mathcal M_1
\otimes \mathcal M_2  \to  \nsets ( \mathcal N_1 \otimes \mathcal N_2) $   by
  $f(\sum v_i\otimes w_i ) = \{ \sum f_1(v_i' )\otimes f_2(  w_i' ):    \sum v_i'\otimes w_i'  = \sum v_i\otimes w_i  \}. $

  It is easy to see that if $f_1,f_2$ are $\supseteq$-morphisms on hyperpairs then $f_1\otimes_\supseteq f_2$ is a $\supseteq$-morphism. The   difficulty with this observation is that for the hyperpairs of greatest interest (for quotient hyperfields) the natural morphisms are $\subseteq$-morphisms.
\end{remark}

   \subsection{Tensor   extensions of $\preceq$-morphisms and weak morphisms  of pre-semiring pairs}$ $

 One can push this theory further,  taking tensor extensions of pre-semirings as in Definition~\ref{tensext}.

\begin{thm}\label{wa}
   Suppose that $\mcA$ is a weakly admissible  $\tT $-semiring, $\tT' $ is a monoid containing  $\tT $, and $f: \mcA_1\to \mcA_2$ is a weak pre-semiring morphism  (resp.~ $\preceq$-morphism).
  If $\tT' $
is the disjoint union  $ \underset{\cdot}{\cup}\,   c_i\tT$ for $c_i\in \tT',$ then    $f$ extends naturally to a
    weak pre-semiring morphism  (resp.~ $\preceq$-morphism) $\tilde f :  {\tT'}\otimes _{\tT} \mathcal M_2 \to   \tT'  \otimes _{\tT} \mathcal N_2 $ given by $\tilde f ( c_i a\otimes y) = c_i\otimes f(ay).$
\end{thm}

  \begin{proof} We do it for $\preceq$-morphisms.

 \begin{equation}
        \begin{aligned}
              \text{(i)} \quad  \tilde f \left(\sum _i( a_i \otimes   b_i ) +\sum _i( a_i' \otimes  b_i )\right)=\tilde f\bigg( \sum _i( a_i + a_{i}') & b_i\bigg) \\ \preceq \sum _i a_i \otimes f(b_i) +\sum _i a_{i}'  \otimes f(b_i)= \tilde f (\sum _i  a_i  & \otimes  b_i ) + \tilde f (\sum _i a_i' \otimes  b_i ), \quad \forall a_i,a_i'\in \tT'.
        \end{aligned}
    \end{equation}

     (ii) Using Theorem~\ref{tenfa}, we have
\begin{equation}
    \begin{aligned}     \qquad
       \tilde f (c_i a \otimes y + & c_i a' \otimes y' ) =    \tilde f (c_i (a \otimes y +   a '\otimes y' )) = c_i \otimes f(ay+a'y') \\& \preceq c_i a \otimes f(y)+c_i a' \otimes f(y')=  \tilde f (c_i a \otimes y) +\tilde f(c_i a' \otimes y' ) .
    \end{aligned}
\end{equation}
  \end{proof}

\begin{example}\label{te}
Any $\preceq$-morphism (resp.~weak morphism)  $f: (\mcA , \mcA_0) \to (\mcA',\mcA'_0) $  of semiring pairs extends to $\preceq$-morphism (resp.~weak morphism) $f:(\mcA [\lambda]  , \mcA_0[\lambda] ) \to (\mcA'[\lambda] ,\mcA'_0[\lambda] ) $ of the polynomial semiring pairs, taking $\tT' = \cup_{i\in \Net} \tT \lambda^i.$
\end{example}

\subsection{The adjoint correspondence}$ $

\begin{rem}\label{adjm}

Suppose that  $\mcM _i$ are $(\tT_i,\tT_{i+1})$-bimodules. The celebrated adjoint isomorphism theorem says  $$\Hom ( \mcM_1 \otimes \mcM_2, \mcM_3) \approx \Hom  (   \mcM_2,  \Hom  ( \mcM_1 ,\mcM_3)).$$
The usual proof for modules over rings, say in \cite[Proposition~2.10.9]{Row06} matches values on simple tensors.  Given $f\in \Hom ( \mcM_1 \otimes \mcM_2, \mcM_3)$ and $w\in \mcM_1$,
one defines $g_w :\mcM_1\to \mcM_3$ by $g_w(w_1) = f(w\otimes w_2).$  In the other direction,
given $g: \mcM_2 \to \Hom  ( \mcM_1 ,\mcM_3) $ one defines the homomorphism $f :\mcM_1 \otimes \mcM_2 \to \mcM_3$ from its action on simple tensors, $f(w_1\otimes w_2) = (g(w_2))(w_1)$. This argument does not use negation, and thus works also for homomorphisms of modules over semirings.
\end{rem}

But what about different sorts of morphisms, i.e., $\preceq$-morphisms or weak morphisms? First of all, we have seen for $\preceq$-morphisms that $f\otimes g$
need not be well-defined on simple tensors. Even if it is well-defined, extending it to all tensors is a challenge, which we take on in this section.
First we settle for a weaker assertion.

 \begin{thm}\label{adj1}
     $\WMor ( \mcM_1 \times \mcM_2, \mcM_3) \approx \WMorp  (   \mcM_2,  \WMorp  ( \mcM_1 ,\mcM_3)).$
 \end{thm}
\begin{proof}
      As in the proof of \cite[Proposition~2.10.9]{Row06}. Given $f\in\WMor ( \mcM_1 \times \mcM_2, \mcM_3)$ and  $w\in \mcM_2,$ we define $f_w : \mcM_1 \to \mcM_3$ by $v \mapsto f(v\otimes w).$ Then for $\sum v_i \in {\mcM_1}_0$,
      $$f_w(v_1+v_2)= f \left(\left(\sum v_i\right)\times w\right)= f\left(\sum v_i  \times w \right) \in {\mcM_3}_0,$$  and likewise for the other side, so $f_w\in  \WMorp  ( \mcM_1 ,\mcM_3).$
    We get $\Phi(f) \in \WMorp  (   \mcM_2,  \WMorp  ( \mcM_1 ,\mcM_3))$ by declaring $\Phi(f)(w)= f_w.$

    Conversely, given $g: \mcM_2 \to \WMorp(\mcM_1,\mcM_3) , $ one   defines the  map $\psi_g:(v,w)\mapsto g(w)(v)$.

    These correspondences clearly are inverses.
\end{proof}

 However,  the maps thus obtained in $\WMor ( \mcM_1 \times \mcM_2, \mcM_3)$
 need not be balanced, so might not produce maps of tensor products.

 The following general observation specializes to the usual adjoint isomorphism in module theory.

\begin{lem}\label{adjc}     $\mcM _i$ are $(\tT_i,\tT_{i+1})$-bimodules. Define $\Morp _ {\tT_1}( \mcM_1 \otimes \mcM_2, \mcM_3)$ to be the  $\preceq$-morphisms. Then there is a reverse map
   $$\Phi :\Morp _ {\tT_1}( \mcM_1 \otimes \mcM_2, \mcM_3) \to \Morp _ {\tT_2}(   \mcM_2,  \Morp _ {\tT_1}( \mcM_1 ,\mcM_3)),$$  given in the proof,
   and, when $\mcM_3$ is a lattice, a map $$\Psi :\Morp _ {\tT_2}(   \mcM_2,  \Morp _ {\tT_1}( \mcM_1 ,\mcM_3)) \to \Morp _ {\tT_1}( \mcM_1 \otimes \mcM_2, \mcM_3).$$
\end{lem}
\begin{proof}
    As in the proof of \cite[Proposition~2.10.9]{Row06}. Given $f\in \Morp _ {\tT_1}( \mcM_1 \otimes \mcM_2, \mcM_3)$ and an  $w\in \mcM_2,$ we define $f_w : \mcM_1 \to \mcM_3$ by $v \mapsto f(v\otimes w).$ Then $$f_w(v_1+v_2)= f ((v_1+v_2)\otimes w)= f(v_1 \otimes w + v_2 \otimes w) \preceq f(v_1\otimes w) + f(v_2\otimes w),$$ so $f_w\in  \Morp _ {\tT_1}( \mcM_1 ,\mcM_3).$
    We get $\Phi(f)\in \Morp _ {\tT_2}(   \mcM_2,  \Morp _ {\tT_1}( \mcM_1 ,\mcM_3))$ by declaring $\Phi(f)(w)= f_w.$

    Conversely, given $g: \mcM_2 \to \Morp(\mcM_1,\mcM_3) , $ one can define the  map $\psi_g:(v,w)\mapsto g(w)(v)$, which yields a right $\preceq$-defined morphism $\bar\psi_g$, namely $$\bar\psi_g (\sum _i v_i\otimes w_i) = \bigwedge \left\{\sum g(w'_i)(v'_i):  \sum  v_i\otimes w_i = \sum  v'_i\otimes w'_i\right\}.$$
\end{proof}

Note that hyperpairs are lattices, where the inf is set intersection.

\section{Directions for further research}\label{FR}

Here are some questions which hopefully could lead to further developments.

\begin{enumerate}
    \item What is the ideal structure of the tensor product of semiring pairs? In particular, of tensor extensions?

    \item Define the \textbf{matrix pair}
$(M_n(\mcA),M_n(\mcA_0))$, where $S$ is the monoid of {\it matrix units} $e_{i,j},$ which have $1$ in the $i,j$ position and $0$ elsewhere. The underlying set of tangible elements of $M_n(\mcA)$ is $\cup_{i,j} \tT e_{i,j}$, cf.~Remark~\ref{cl0}. Is the tensor product of matrix pairs a matrix pair?

      \item (Suggested by the referee.) There is a well-known process of tropicalization  of fields  described in \S\ref{ME}(vii).  How does the tensor product of tropicalizations  compare to the tropicalization of the tensor product? Again, this seems more accessible if one considers tropicalizations of tensor extensions.
\end{enumerate}

\end{document}